\documentclass[14pt]{extarticle}
\usepackage[margin=1in]{geometry}
\usepackage{dsfont}
\usepackage{amsfonts}
\newtheorem{theorem}{Theorem}[section]
\newtheorem{pro}{Proposition}[section]
\newtheorem{lemma}{Lemma}[section]

\newtheorem{cor}{Corollary}[section]

\newcommand{\proof}[1]{\noindent{\it\bf Proof:#1\ }}
\newcommand{\QED}{\hfill$\Box$\medskip}

\begin{document}

\title{On    Linear Part of Filled-Section in Splicing }
\author{
\\
Gang Liu \\Department of Mathematics\\ UCLA}
\date{July 27,  2018}
\maketitle

\section{ Introduction}
In \cite{1} and \cite{2}, using splicing  Hofer, Wysocki and Zehnder have introduced a smooth model for GW theory   near  a  nodal  map $f:S=C_-\cup C_+\rightarrow M$
in the setting of Sc-Banach manifolds in  polyfold theory. The basic idea is to introduce, for each gluing parameter $R=:R_{\theta}\in [R_0, \infty)$ with $\theta\in S^1$,  the anti-gluing $S_-^R=:C_-\ominus_R C_+$ as  the counter part of the usual  gluing
$S^R_+=:C_-\oplus_R C_+$; then to define a family of isomorphisms $T^R$/${\hat T}^R$ between the spaces of  maps/sections over $S$ to the corresponding ones over $S^R=S^R_-\cup S^R_+.$ The usual  ${\overline {\partial}}_J$-operator, $\Phi^R_+$, acting on the maps  on $S^R_+$ then  is  extended into an operator, denoted by $\Phi^R=(\Phi^R_-, \Phi^R_+)$ acting on the maps  on $S^R$ in such a way that  the index of the linearizations   of $\Phi^R_-$  is equal to zero so that  the zero set of $\Phi^R$, at least virtually, can be identified with the usual moduli space of $J$-holomorphic maps in the GW-theory  with the domain $S^R_+$.
The filled-section $\Psi=\{\Psi^R\}$
is obtained from $\Phi^R$ by the conjugation:
$\Psi={\hat T}^{-1}\circ \Phi\circ T=\{{\hat T}^{-1}\circ \Phi^R\circ T\}$.
Unlike the ${\overline {\partial}}_J$-operators $\{\Phi^R_+\},$
the domains and targets of  $\{\Psi^R\}$ are  independent  of
$R$, this make it possible to take derivatives for these filled-sections along $R$-direction even at $R=\infty$.  However since  $T$/${\hat T}$ is only an Sc-isomorphism, it is expected that there is a loss of differentiability for  $\Psi$ in general. In actuality,  the derivatives of filled-section defined in \cite{1,2}
have  a loss of differentiability in the  nonlinear part of $\Psi$, and do not 
converge fast enough as $R\rightarrow \infty$ in the  ``linear part''.  
This together with other difficulties prevents from proving that $\Psi$ is smooth in the usual setting of Banach analysis. Instead   it  was  proved in \cite{ 2} that 
the filled -section there   is Sc-smooth  with respect to  the  gluing profile $R=\phi(r)=e^{1/r}-e^{1/r_0}$.

Since in the splicing above the operator  $\Phi^R_-$  is free to be chosen as long as some basic requirements are  satisfied,
 %At this  point, we note that
% Notice that  while  the operator $\Phi^R_+$ acting on the finite cylinder $C_+^R$ is required to be the standard ${\overline {\partial }}_J$-operator,
it is natural to  ask what  the most natural choices are
and if
 it is possible to define the filled-section using such choices  in the framework of the usual Banach analysis.
 This sequence of papers is a report of our work with an affirmative to this question.

  While  the operator $\Phi^R_+$ acting on the maps on finite cylinder $S_+^R$ is required to be the standard ${\overline {\partial }}_J$-operator,      the desired
    operator  $\Phi^R_-$  on $S^R_-\simeq {\bf R}^1\times S^1$ cannot act on  $u_-\ominus u_+$ alone
    and has to be nonlinear (comparing with  \cite{ 2}).   In cylindrical coordinated $(t, s)$ on $S^R_-$,
$\Phi^R_-$ will have  the form $\Phi^R_- (u_-\ominus_R u_+)=\partial_t (u_-\ominus_R u_+)+ \omega^R_-(t) (u_-\ominus_R u_+, u_-\oplus_R u_+)+J(u_-{\hat \oplus}_R u_+\circ \Gamma)\partial_s(u_-\ominus_R u_+).$ Here $\Gamma^R:S^R_-\rightarrow {\hat S}^R_+$ is the {\bf transfer map} that transforms the almost  complex structure along $u_-{\hat \oplus}_R u_+$ to the one along $u_-\ominus_R u_+$, and $u_-{\hat \oplus}_R u_+$ is the extended gluing of $u_-$ and $u_+$ with domain ${\bf R}^1\times S^1$; and $\omega^R_-(t)$ is part of the connection matrix $\omega^R(t)$ with $\omega^R_+(t)=0$ so that
 the linear part $\partial_t$ of $\Phi^R $  from the trivial connection  is replaced by
  the covariant derivative from the $(R, t)$-dependent connection given by $\omega^R(t)$
  such that on  $u_-\oplus_R u_+$, it is still same as $\partial_t$ (see the definition in next section).
  
  The ideas to deal with  linear and nonlinear part of $\Phi^R$ are quite different. In this paper we only   define the linear part of $\Phi^R$,
  denoted by $\Phi_L^{R}=((\Phi_L^{R})_-, (\Phi_L^{R})_+)$.
  Since $ (\Phi_L^{R})_+=\partial_t$,   its study has its own interests  in analysis beyond the need in the Gromov-Witten
  type theories.
  
  Let $\Psi^{R}_L=(T^R)^{-1}\circ \Phi_ L^{R}\circ T^R$  be the corresponding operator.
  Then  the main theorem of this paper is the following theorem.

  \begin{theorem}
  	Using the gluing profile $R=e^{1/r}-e^{1/r_0}$, the filled-section $\Psi_L=\{\Psi^{R_{\theta}}_L\}$ above  with $r\in [0, r_0)$ (hence $R=R_{\theta}\in (R_0, \infty]$ ) and $\theta\in S^1$ is of class $C^1$.
  \end{theorem}

  We now briefly explain the main idea of the definition  of  $\Psi^{R}_{L}$. The main  issue here is to  make the  proper choice of $\Phi^R_{L, -}$ so that the  derivatives of $\Psi^{R}_{L}$ converge or decay fast enough.

  To this end,  we need change the matrix of splicing $T_{\beta}$ used to define the
   total gluing map $T^R$: instead of using cut-off functions $\beta$ and $1-\beta$ as entries of $T_{\beta}$  to define $T^R$ based on partition of unit (at least for $u_-\oplus u_+$) in  \cite{1,2} with the support of ${\beta}'$ lying on the sub-cylinder of length $1$ with the distance  $R$ to the boundary circles of $C_{\pm}$,
    a pair of $R$-dependent cut-off functions $\beta^R=(\beta^R_-, \beta^R_+)$  is chosen  without the condition that $\beta^R_-+\beta^R_+
   =1$ in such a way that (i) $|(\beta^R_{\pm})'|\sim o(1/(R^{1/2}\cdot \ln^2 R)$; (ii)$(\beta^R_{\pm})'$   has the support  with length $\sim 2R^{1/2}\cdot \ln^2$  with  the distance $\sim (R-3R^{1/2}\cdot \ln^2)$ to the nearest boundary circle of $S^R$   so that the splicing matrix $T_{\beta}$ has two
  ``splicing regions'',  where no-trivial splicing takes place (see the definition in next section).
   
    The purpose of this construction of the splicing matrix is to create a situation such that
   (i) and (ii) above together has the following effect on the "unwanted terms" in the  derivatives of $\Psi^R_{-, L}$:  they become favorable terms in one splicing region  but have  worse convergent rate in the other. The key  step then is  to  introduce  the connection matrix $\omega^R$ to get rid of  these new unwanted terms.

%In particular, the splicing matrix $T_{\beta}$ has two
% regions where no-trivial splicing takes place (see the definition in next section)  and the  curvatures of the connection matrices $\omega^R(t)$ are not equal to zero but localized on those places
% ( comparing with  \cite{ 2} ). 

We remark that  using the splicing matrix $T_{\beta}$ and the connection matrices $\omega^R(t)$ above, 
 the linear splicing
 here   permits  the  interpretation  as a pair of traveling  rank 2 (or complex line)
bundles with non-trivial curvatures (see Sec. 2.5). The more details of this and  its relation with the corresponding linear splicing in the setting of \cite{2} will be given  somewhere else.

 Theorem 1.1  will be proved in Sec. 3 and Sec. 4 after the basic definitions  of splicing are given  in Sec.2.

The nonlinear part of the filled-section in the setting of Banach analysis  will
be given in \cite{3}.

Only elementary  Banach analysis is used in this paper which can be found in Lang's book \cite{4}. In last section, we have collected some  of  basic analytic  facts used in this paper and its companion \cite{3}.

\section{Basic definitions of the  splicing}

The basic analytic  set-up  used in  this paper and its companion is the space 
$L_{k, \delta}^p(C_{\pm }, {\bf C}^n)$ of $L_{k}^p$- functions with  $\delta$-exponential decay.   We use such a space  as a local model  for  $L_k^p$-maps near a nodal stable map with fixed values at its ends.
In order to have the ``right'' analytic set-up to accommodate the geometric operator
$\Phi$, in particular to have the right dimension of the zero set of $\Phi^R_+$,
it is necessary  to allow the ends moving; 
in the case with fixed ends  the constrains should be imposed to the image of the total gluing  accordingly. However, in this paper and its companion, 
we  will suppress   this and related  aspects in order to concentrated on  the main  issues here: defining  the filled-section, for which the   total gluing, the main construction  of this section  only  serves   as an intermediate step. 

 %The general case allowing the ends moving will be treated somewhere else.  Our work here then established the analysis along each ``fiber'' for the general case.

 The constructions in subsections 2.1 and 2.2 are essentially  the same as    those  in \cite{2}.
  %The transfer map $\Gamma^R$ is introduced in subsection 2.3 and
 % The linear part ${\Psi_L}$ of the  filled-section is defined in 2.4.
\subsection{Total gluing of the nodal surface  $S$}
We start with  the  definitions of the total gluing of the given   nodal surface $S$.

%In the following, we will work with the usual logarithmic gluing profile $t=-\ln |z|$ of ${\overline {\cal M}}_{g, k}$.
%Our results will be formulated in term of variable $t>0$ and  asymptotic analysis as $t$ goes to infinity.

For the purpose of  this paper, we only need to consider  the  germ   of  the given   nodal  surface, still denoted by 
  $S=C_-\cup_{d_-=d_+}C_+$   with the  double point
 ${d_-=d_+}$. Thus
 each  component $(C_{\pm}, d_{\pm})$ is identified with the standard disk  $ (D, 0)$
 canonically upto a rotation.  Identify $(C_{-}, d_{-})$ with $ ((-\infty, 0)\times S^1, -\infty\times S^1 )=(L_-\times S^1, -\infty\times S^1)$
 canonically upto a rotation by  considering the double point $d_-$ as the $S^1$ at $-\infty$ of the half cylinder $L_-\times S^1$. Here we have denoted the negative half line $(-\infty, 0)$ by $L_-$.
 Similarly $(C_{+}, d_{+})\simeq (( 0, \infty,)\times S^1, \infty\times S^1 )=(L_+\times S^1, \infty\times S^1).$

 \medskip
 \noindent ${\bf \bullet}$  {\bf Cylindrical coordinates on $C_{\pm}:$}

\medskip
\noindent
  By the identification  $C_{\pm}\simeq L_{\pm}\times S^1$,   each $C_{\pm}$ has  the  cylindrical coordinates $(t_{\pm}, s_{\pm})\in L_{\pm}\times S^1.$

Let $a=(R, \theta)\in [0, \infty]\times S^1$ be the gluing parameter. 
 To defined the total gluing/deformation $S^a=S^{(R, \theta)}$ with gluing parameter $R\not= \infty$, we  introduce  the $a$-dependent cylindrical coordinates $(t^{a, {\pm}}, s^{a,{\pm}})$ on  $C_{\pm}$ by the formula 
  $t_{\pm}=t^{a, \pm }\pm  R$ and $s_{\pm }=s^{a, \pm}\pm\theta$. In the following if there is no confusion, we will  denote  $t^{a, \pm }$ by $t$ and $s^{a, \pm}$ by $s$ for both of these $a$ -dependent cylindrical coordinates.

 Thus the $t$-range for $L_-$ is $(-\infty, R)$ and  the $t$-range for $L_+$ is $(-R, \infty)$ with the intersection $L_-\cap L_+=(-R, R)$.

 \medskip
 \noindent ${\bf \bullet}$  {\bf Total gluing $S^a=(S^a_-, S^a_+)$ :}

 \medskip
 \noindent

%The construction here depends on a choice of a length function $L$ that is fixed.
  In term of the $a$-dependent cylindrical coordinates $(t,s)$,   $C_-=(-\infty, R)\times S^1$ and $C_+=(-R,\infty)\times S^1$.
  
  Then  $S^a_+$ is defined to be the finite cylinder of length $2R$  obtained by gluing $(-1, R)\times S^1\subset C_-$
with $(-R, 1)\times S^1\subset C_+$ along the "common" region $(-1, 1)\times S^1$ by the identity map in term of the $a$-dependent  coordinates $(t,s)$.
 Similarly,  $S^a_-$ is the infinite cylinder defined by  gluing $(-\infty, 1)\times S^1\subset C_-$
with $(-1, \infty)\times S^1\subset C_+$ along  $(-1, 1)\times S^1$ by the identity map.

Geometrically, both $S^a_{\pm}$ are obtained by first cutting  each $C_{\pm}$ along the circle at $t=0$ into two sub-cylinders,then gluing back  the  sub-cylinders in $C_-$ with  the corresponding ones in $C_+$ along  the same circle  with a relative rotation of angle $2\theta$.
Set $S^{\infty}=S.$

The total gluing can also be described as follows. We may consider $S=S^{\infty}=C_-\cup C_+=L_-\times S^1\cup L_+\times S^1$  to be the surface by gluing $C_-$ and $C_+$ along the the circles $\{{-\infty}\}\times S^1$ and  $\{{\infty}\}\times S^1$ by identity map. Then
$S^{\infty, \theta}$ is defined to be the surface still  by  gluing $C_-$ and $C_+$ along the these  circles but by a relative rotation of angle $2\theta.$ Then $S^a$ is defined to be $S^a=(S^{\infty, \theta})^R.$

 Now the  cylindrical coordinates $(t_{\pm}, s_{\pm})$ on $C_{\pm}$ as well as the $a$-dependent cylindrical coordinates
 $(t, s)$  become the corresponding ones on each  $S^{a}_\pm$ with the relation:   $t_\pm=t\pm R$ and $s_\pm =s\pm \theta.$

\subsection{ Splicing matrix by HWZ }

 The splicing matrix  used in \cite{1, 2}   will be denoted by $T_{\alpha}$  which is defined  as follows.

Fix a smooth  cut-off function $\alpha:{\bf R}^1\rightarrow [0, 1]$ with the property
that  $\alpha  (t)=1$ for $ t<-1$,  $\alpha (t)=0 $ if $ t>1 $ and $\alpha'\leq 0$.

Then $$T_{\alpha}=\left[\begin{array}{ll}
\alpha &  -(1-\alpha)\\
(1-\alpha)   & \alpha
\end{array}\right].$$

Though it is not necessary, $\alpha$ can be chosen  such  that
$\alpha(t)+\alpha(-t)=1.$

%$$ \alpha (t)= \left\{ \begin{array}{ll}
%\beta  & 1-\beta \\
%-(1-\beta ) & \beta
%
\subsection{ Splicing matrix $T_{\beta}$}
To defined   $T_{\beta}$,   we need to choose a
length function depending $R$. For the purpose of this paper, the length fuction  $L(R)=L_1(R)=R^{1/2} \cdot ln^2 R$. In general, for any positive integer $k$, $L_k(R)=R^{k/(k+1)}\cdot ln^2 R.$

The  splicing matrix $T_{\beta}$ used in this paper is defined by using a new pair of  cut-off function $\beta=(\beta_-,\beta_+)$ depending on the two  parameters $(l, d)$ that parametrize the group of affine transformations $\{t\rightarrow lt+d\}$ of ${\bf R}^1$ defined as follows.

 Rename $\alpha$ by $\alpha_-$ and $1-\alpha$ by $\alpha_+$. Fix  $l_0>1$ and $d_0>1$. Then  $\beta_{\pm}=\{\beta_{\pm;l,d}\}:{\bf R}^1\times [l_0, \infty)\times [d_0, \infty)\rightarrow [0, 1]$   defined by $\beta_{\pm }(t, l,d)=\alpha_{\pm}( (t\pm d)/l)$, or $\beta_{\pm ;l,d}=\rho_l\circ \tau_{\pm d}\alpha$. Here the translation and multiplication operators are  defined   by $\tau _{d} (\xi)(t)=\xi(t+d)$ and $\rho_l(\xi)(t)=\xi(t/l)$ respectively.

The pair $(l, d)$ will  be  the functions on $R$, $(l=L(R),d=d(R)) $ defined above with   $d\geq 3l$.  To be specific,  set $d=3l.$
 Clearly $\beta_{\pm}$ is a smooth cut-off function  with the following two properties:

 \medskip
 \noindent
 $P_1:$ for  $k\leq k_0$ the $C^0$-norm of the $k$-th derivative $\|\beta_{\pm}^{(k)}\|_{C^0}\leq C/l^k$  , where $C=\|\alpha\|_{C^{k_0}},$
  which is independent of $l$;

 \medskip
 \noindent
 $P_2:$ under the assumption that $d\geq 3l$, the support of $\beta_{-}'$  is contained in the interval$(d-l,  d+l)$ with  $\beta_{-}=1 $ on   $(-\infty,d-l]$ and $\beta_{-}=0 $ on   $[d+l, \infty);$ and $\beta_{+}'$  is contained in the interval$(-d-l,  -d+l)$ with  $\beta_{+}=1 $ on   $[-d+l, \infty)$ and $\beta_{+}=0 $ on   $(-\infty, -d-l]$.

 Thus comparing with $\alpha_-$ and $\alpha_+$,  whose  supports are in the common interval $[-1, 1]$,  the supports of  $\beta_-'$ and $\beta_+'$ are in the two intervals $[-d-l, -d+l]$ and $[d-l, d+l]$ without overlaps.
 Those intervals  are corresponding to the regions where the  splicing takes place.
 %Thus  by using the standard matrix $T_{\alpha}$  is  a quite different from the one using the non-standard $T_{\beta}$ 
 The splicing matrix then is  defined by

 $$T_{\beta}=\left[\begin{array}{ll}
 \beta_- &  -\beta_+\\
 \beta_+   & \beta_-
 \end{array}\right]. $$

 Note that from $P_2$, on $(-d+l, d-l)$, $\beta_-=\beta_+=1$.
 Then  for $t$  in the  three intervals
 $(-\infty, -d-l)$,$(-d+l, d-l)$ and $(d+l, \infty)$,  $T_{\beta}(t)$    are the following constant matrices

 $$M_1=Id=\left[\begin{array}{ll}
 1 &  0\\
0   &1
 \end{array}\right], M_2=\left[\begin{array}{ll}
 1&  -1\\
 1  & 1
 \end{array}\right], \, \, and\,\, M_3= \left[\begin{array}{ll}
 0 & -1\\
 	1  & 0
 \end{array}\right] .$$

 The length $2(l+d)$ of the interval $(-l-d, l+d)$  with $l=L(R)$ is defined to be the {\bf length of the splicing  matrix $T_{\beta}$}.

 % Note that  $T_{\alpha}$ only has   two constant matrix $M_1$ and $M_3$ above outside the interval $ (-1, 1)$; the  constant matrix $T_{\beta}=M_2$  on the middle interval
 %$(-d+l, d-l)$ gives rise  a total gluing map $T^R_{\beta}$  in the region  which has no comparison with $T^R_{\alpha}$.

  Note that  $\beta_{\pm}(t)<1$ implies that $\beta_{\mp}(t)=1$ so that

the determinant   $$1\leq D= det\, \left[\begin{array}{ll}
\beta_-  & -\beta_+ \\
\beta_+ & \beta_-
\end{array}\right]   =\beta_-^2+\beta_+^2\leq 2.$$

This implies that $T^a= (\ominus_a, \oplus_a)$ defined below is invertible uniformly.

\subsection{Total gluing  $T^a$ of maps and sections}

%Since the matrix $T_{\beta}$ for the total gluing $T^R=(T^R_-, T^R_+)=(\ominus_R,\oplus_R )$ is $s$-independent, we  define the total gluing $T^R$ on the maps and sections on $C_{\pm}$ first.

Let $C^{\infty }_c(C_{\pm}, E )$ be the set of $E$-valued $C^{\infty}$ functions 
on $C_{\pm}$ with compact support, where $E={\bf C} ^n;$  similarly for $C^{\infty }_c(S^a_{-}, E).$
%consists of all  $E$-valued smooth functions on $S^a_{\pm}$ with compact support. 
Note that for the definitions here  the surfaces $C_\pm$ and $S^a_{-}$  are considered as the cylinders with boundary.

Then $T^a=(T^a_-, T^a_+):  C^{\infty}_c(C_{-} E)\times C^{\infty} _c(C_{+}, E)\rightarrow C^{\infty}_c(S^a_{-} E)\times C^{\infty}(S^a_{+}, E)$ is defined as follows.

In matrix notation, for each $(\xi_-, \xi_+)\in  C^{\infty}(C_{-}, E)\times C^{\infty}(C_{+}, E)$ considered as a column vector,
$$ T^a((\xi_{-}, \xi_{+}))=(T^a_-(\xi_{-}, \xi_{+}), T^a_+(\xi_{-}, \xi_{+}))=(\xi_{-}\ominus_a\xi_{+},  \xi_{-}\oplus \oplus_a \xi_{+})$$ $$ =\left[\begin{array}{ll}
\beta_- &  -\beta_+\\
\beta_+   & \beta_-
\end{array}\right]
\left[\begin{array}{c}
\tau_{-a}\xi_{-} \\
\tau _{a} \xi_{+}
\end{array}\right]
.$$

 The inverse of the total gluing,  $(T^a)^{-1}=(T^a_-, T^a_+)^{-1}:C^{\infty}_c(S^a_{-}, E)\times C^{\infty}(S^a_{+}, E)  \rightarrow C^{\infty}_c(C_{-}, E)\times C^{\infty}_c(C_{+}, E)$ is defined as following:
for a pair of the $E$-valued  functions $(\eta_{-}, \eta_{+})\in C^{\infty}(S^a_{-} E)\times C^{\infty}(S^a_{+}, E)$,

$(T^a)^{-1}(\eta_{-}, \eta_{+})=(\oplus_a \oplus \ominus_a)^{-1}(\eta_{-}, \eta_{+})=$
$$\left[\begin{array}{ll}
\tau_{a} & 0 \\
0  &  \tau_{-a}
\end{array}\right] \cdot{\frac {1}{D}}\cdot
\left[\begin{array}{ll}
\beta_- &  \beta_+\\
-\beta_+  & \beta_-
\end{array}\right]
\left[\begin{array}{c}
\eta_{-} \\
\eta _{+}

\end{array}\right].$$

%For $a=Re^{i\theta}=(R, \theta)$, the total gluing map
% $T^a=(T^a_-, T^a_+): C_c(C_{-},  E)\times C_c(C_{+}, E) =: C_c(L_{-}\times S^1,  E)\times C_c(L_{+}\times S^1, E)\rightarrow C_c(C^a_{-}, E) \times C_c(C^a_{+}, E)  $ is defined by the formula
% above but replacing $\tau_R$ by $\tau_a=\tau_R\circ \tau_{\theta}$.

%\medskip
%${\bullet}$ Note on notations:

 %In the following we  also write  $a$ as  $R_{\theta}$.
 %In this notation,
%$\tau_{-R_{\theta}}=\tau_{-R}\circ \tau_{-\theta}=(\tau_{a})^{-1}.$

%Note that the map $(u_-, u_+)\in C^{\infty}(L_{-}\times S^1,  E)\times C^{\infty}(L_{+}\times S^1, E)$ satisfies the asymptotic condition that
%$u_-(-\infty)=u_+(\infty)$ that is corresponding to the condition that $T^a_-(u_-, u_+)(-\infty)=-T^a_-(u_-, u_+)(+\infty)$.
%In other words $T^a$ maps these two subspaces each other isomorphically as before.

\subsection{Complex total gluing  }

Denote $\tau_{-a}\xi_{\pm}$ in the previous subsection  by $\eta_{\pm}$.
Then the above formula for the total gluing can be express in term of complex
 notations as follows.
 
 Let $\eta=\eta_{-}+\eta_{+}i$ and $\beta=\beta_-+\beta_+i$. Then 
$\beta\cdot \eta=(\beta_-\eta_{-}-\beta_+\eta_{+})+(\beta_-\eta_{+}+\beta_+\eta_{-})i,$

which is the same as 
$$\left[\begin{array}{ll}
\beta_- &  -\beta_+\\
\beta_+   & \beta_-
\end{array}\right]
\left[\begin{array}{c}
\eta_{-} \\
\eta_{+}
\end{array}\right]
$$ by matrix notations.

 \section{ The linear part   $\Psi_L$ of the filled-section}

%\medskip
%\noindent ${\bf \bullet}$   {\bf Length function $L(R) $}:

%\medskip
%\noindent
%The Length function
%$L$ is defined to be  $L(R)=L_1(R):R^{1/2}(ln R)^2$. In general, for a positive integer $k, $  $ L(R)=L_k(R):R^{k/(k+1)}(ln R)^2$.   We will show that   with the choice of this $L_1$, the resulting $\Psi$ is of class $C^1$.

   \medskip
   \noindent
   {\bf The linear part $\Psi^R_L$ and $\Psi^a_L$ :}

  Let $T_{\beta}$  be the matrix associated to the total gluing map $T^a$
   and $$\omega (t)=-\left[\begin{array}{ll}
   \beta'_{-} & -\beta'_{+} \\
   0  &  0
   \end{array}\right] \cdot T_{\beta}^{-1}=-1/D\left[\begin{array}{ll}
   \beta'_{-}\beta_-+ \beta'_{+}\beta_+& \beta'_{-}\beta_+- \beta_+'\beta_-\\
   0  &  0

   \end{array}\right] $$ denoted by $$\left[\begin{array}{ll}
   	e& f\\
   	0  &  0
   \end{array}\right] .$$

    Note that in above, we have suppressed the $a$-dependence of $e$, $f$ and $\omega$ in our notations, and we will continue to do so.
   Then $$\Phi_L^a\left[\begin{array}{c}
    {u}_{-}\ominus u_+  \\
   {u}_{-}\oplus u_+
   \end{array}\right] =( {\partial_t +\omega (t)}) \left[\begin{array}{c}
   {u}_{-}\ominus u_+  \\
   {u}_{-}\oplus u_+
   \end{array}\right]=\partial_t\left[\begin{array}{c}
   {u}_{-}\ominus u_+  \\
   {u}_{-}\oplus u_+
   \end{array}\right]
     +\omega (t) \left[\begin{array}{c}
   {u}_{-}\ominus u_+  \\
   {u}_{-}\oplus u_+
   \end{array}\right]$$
   $$=\partial_t \left[\begin{array}{c}
   {u}_{-}\ominus u_+  \\
   {u}_{-}\oplus u_+
   \end{array}\right]+\left[\begin{array}{c}
   e{u}_{-}\ominus u_+ +f{u}_{-}\oplus u_+ \\
    0
   \end{array}\right].$$

    This prove the following
  \begin{lemma}
  For $R\not=\infty$, $(\Phi_L^a)_+({u}_{-}\oplus u_+ )=\partial_t ({u}_{-}\oplus u_+)$   as required, and 	 $(\Phi_L^a)_-({u}_{-}\ominus  u_+ )=\partial_t ({u}_{-}\ominus u_+)+ \omega(t) ({u}_{-}\ominus  u_+,{u}_{-}\oplus u_+ ).$
  \end{lemma}

\begin{lemma}
 $$\Psi^a_L(u_-, u_+)={\partial_t }\left[\begin{array}{c}
  {u}_{-} \\
  u _{+}
  \end{array}\right] +
  \left[\begin{array}{c}
 \tau_{a}( \beta_+\beta_+'/D) u_-   + \tau_{a} (\beta_+\beta_-'/D){\tau_{2a}	u_{+}}\\
\tau_{-a} ( \beta_-\beta_+'/D)\tau_{-2a} {u}_{-} + \tau_{-a}(\beta_- \beta_-'/D)u_+ \end{array}\right]  .$$
 
 There is no loss of  differentiability in $\Psi_L^a(u_-,u_+)
  $. For $t\not \in (-d-l+1, d+l-1)$,
  $\Psi^a_L(u_-,u_+)=\partial _t(u_-,u_+)$.
 \end{lemma}
 
   \proof
   
    By definition $\Psi^a_L=(T^a)^{-1}\circ \Phi_L^a \circ T^a $, and
  $$\Psi^a_L(u_-, u_+)=\left[\begin{array}{ll}
  \tau_{a} & 0 \\
  0  &  \tau_{-a}
  \end{array}\right]
  \cdot T_{\beta}^{-1}\cdot  ( {\partial_t +\omega (t)}) \left \{T_{\beta}\cdot
  \left[\begin{array}{c}
  \tau_{-a} {u}_{-} \\
  {\tau_{a}	u _{+}}
  \end{array}\right] \right\} $$

  $$=\left[\begin{array}{ll}
  \tau_{a} & 0 \\
  0  &  \tau_{-a}
  \end{array}\right]
  \cdot T_{\beta}^{-1}\cdot  \left \{ {\partial_t }T_{\beta}\cdot
  \left[\begin{array}{c}
  \tau_{-a} {u}_{-} \\
  {\tau_{a}	u _{+}}
  \end{array}\right]+
  \omega (t)\cdot   T_{\beta}\cdot
   \left[\begin{array}{c}
   \tau_{-a}  {u}_{-} \\
   {\tau_{a}	u _{+}}
    \end{array}\right]
+T_{\beta}\cdot  \left[\begin{array}{c}
  \tau_{-a}  {\partial_t }{u}_{-} \\
  {\tau_{a} {\partial_t }	u _{+}}
 \end{array}\right]   \right\}
 $$

  $$=\left[\begin{array}{ll}
  \tau_{a} & 0 \\
  0  &  \tau_{-a}
  \end{array}\right]
  \cdot T_{\beta}^{-1}\cdot
  $$ $$
  \left \{ \left[\begin{array}{ll}
  \beta_-' & -\beta_+' \\
  \beta_+' &  \beta_-'
  \end{array}\right] \cdot
  \left[\begin{array}{c}
  \tau_{-a} {u}_{-} \\
  {\tau_{a}	u _{+}}
  \end{array}\right] -\left[\begin{array}{ll}
  \beta'_{-} & -\beta'_{+} \\
  0  &  0
  \end{array}\right]
  \cdot
  \left[\begin{array}{c}
  \tau_{-a}  {u}_{-} \\
  {\tau_{a}	u _{+}}
  \end{array}\right]
   + T_{\beta}\cdot
  \left[\begin{array}{c}
  \tau_{-a}  {\partial_t }{u}_{-} \\
  {\tau_{a} {\partial_t }	u _{+}}
  \end{array}\right]  \right\}$$

  $$={\partial_t }\left[\begin{array}{c}
 {u}_{-} \\
  	u _{+}
  \end{array}\right] +
   \left[\begin{array}{ll}
  \tau_{a} & 0 \\
  0  &  \tau_{-a}
  \end{array}\right]
  \cdot T_{\beta}^{-1}\cdot    \left[\begin{array}{ll}
 0 &0\\
  \beta_+' &  \beta_-' \end{array}\right]  \cdot
  \left[\begin{array}{c}
  \tau_{-a} {u}_{-} \\
  {\tau_{a}	u _{+}}
  \end{array}\right]   $$

$$={\partial_t }\left[\begin{array}{c}
{u}_{-} \\
u _{+}
\end{array}\right] +
\left[\begin{array}{ll}
\tau_{a} & 0 \\
0  &  \tau_{-a}
\end{array}\right]
\cdot   1/D \left[\begin{array}{ll}
 \beta_-&\beta_+\\
-\beta_+ &  \beta_- \end{array}\right]\cdot    \left[\begin{array}{ll}
0 &0\\
\beta_+' &  \beta_-' \end{array}\right]  \cdot
\left[\begin{array}{c}
\tau_{-a} {u}_{-} \\
{\tau_{a}	u _{+}}
\end{array}\right]   $$

 $$={\partial_t }\left[\begin{array}{c}
 {u}_{-} \\
 u _{+}
 \end{array}\right] +
 \left[\begin{array}{ll}
 \tau_{a} & 0 \\
 0  &  \tau_{-a}
 \end{array}\right]
 \cdot 1/D   \left[\begin{array}{ll}
 \beta_+\beta_+' &\beta_+\beta_-'\\
 \beta_-\beta_+' & \beta_- \beta_-' \end{array}\right]  \cdot
 \left[\begin{array}{c}
 \tau_{-a} {u}_{-} \\
 {\tau_{a}	u _{+}}
 \end{array}\right]   $$

  $$={\partial_t }\left[\begin{array}{c}
  {u}_{-} \\
  u _{+}
  \end{array}\right] +
  \left[\begin{array}{c}
 \tau_{a}( \beta_+\beta_+'/D) u_-   + \tau_{a} (\beta_+\beta_-'/D){\tau_{2a}	u_{+}}\\
\tau_{-a} ( \beta_-\beta_+'/D)\tau_{-2a} {u}_{-} + \tau_{-a}(\beta_- \beta_-'/D)u_+ \end{array}\right]  . $$

 For $t\not \in (-d-l+1, d+l-1)$,   ${\beta}_{\pm}$   is independent of $t$ so that $ {\beta}_{\pm}'=0$. 
 
 \QED

 %This proves the following
 Since the splicing matrix $\beta$ is $s$-independent and  the translation operator $\tau_{\theta}$ appeared in the total gluing map $T^a$ commutes with both ${\partial}_t$ and ${\partial }_s$,
 $\tau_{\theta}$ does not affect  analysis here in any essential way. In  the most part of the rest of this section we will only give the details for  the results    using $T^R$ and state the corresponding ones using $T^a.$

Denote  the error term $$ \left[\begin{array}{c}
  	\tau_{R}( \beta_+\beta_+'/D) u_-   + \tau_{R} (\beta_+\beta_-'/D){\tau_{2R}	u_{+}}\\

  	\tau_{-R} ( \beta_-\beta_+'/D)\tau_{-2R} {u}_{-} +
  	 \tau_{-R}(\beta_- \beta_-'/D)u_+ \end{array}\right]$$
  by $E^R(u_-, u_+)=(E^R_-(u_-, u_+), E_+^R(u_-, u_+)).$

  Then $\Psi^R_L(u_-, u_+)=\partial_t (u_-, u_+)+E(u_-, u_+)$.

\medskip
\noindent
$\bullet$  Exponentially weighted $L_k^p$-maps/sections

Starting from next lemma, we will consider the spaces of $L_{k,\delta}^p$-maps/sections. Throughout this paper, we assume that $k\geq 1$ and $p>2.$

 The  spaces of such maps used in this paper  are  defined as follows.

$$L_{k, \delta}^p(C_{\pm}, E)=\{u_{\pm}:C_{\pm}\simeq (\pm R_0, \pm \infty)\times S^1 \rightarrow E\, |\, \|u_{\pm}||_{k,p, \delta}<\infty\}.$$ Here $\|u_{\pm}||_{k,p, \delta}=\|e_{\pm}\cdot u_{\pm}||_{k,p}$ and the weight function $e_{\pm}(t,s)=e^{\delta |t|}$ where $0<\delta <1$ and $1$ (or $2\pi$ depending on how to parametrize $S^1$) is the smallest positive eigenvalues of the self-adjoint operator $i\partial_s $ acting on complex valued functions on $S^1$.

Thus for each $u_{\pm}\in L_{k, \delta}^p(C_{\pm}, E)$ by Sobolev embedding theorem applying to sub-cylinders $[\pm R, \pm R\pm1]\times S^1$, there exists a constant  $C$ such that $|u_{+}(\pm t,s)|\leq C e^{-\delta|\pm t|} $. In particular $\lim_{\pm t\rightarrow \pm\infty} u_{+}(\pm t,s)=A_0=0$. In other words, we only consider the case of $L_{k, \delta}^p$-maps with fixed end (= $0$ of $E$) in this paper.

%This is sufficient as the local isometries   of $E$ at $A_0=0$ induce the identifications of any other such spaces with asymptotic $A_1\in E$ with the one above preserving all analytic structures near the ends.

  \begin{lemma}
   Considered  as  a map $E=\{E^R, R\in [R^0, \infty]\}:L_{k, \delta}^p(S, E)\times [R^0, \infty]\rightarrow L_{k-1, \delta}^p(S, E)$, hence a map $L_{k, \delta}^p(S, E)\times [0, r_0]\rightarrow L_{k-1, \delta}^p(S, E)$, $E$ is continuous.	So is $\Psi_L:L_{k, \delta}^p(S, E)\times [0, r_0]\rightarrow L_{k-1, \delta}^p(S, E).$
  \end{lemma}

  \proof

  For $R\not =\infty$, this is reduced to show  that $F:L_{k, \delta}^p(S, E)\times [0, r_0]\rightarrow L_{k-1, \delta}^p(S, E)$ give by $F(u, R)=u\circ \tau_R$
is continuous. The result is well-known. A proof of this is given in the last section.

Note that for $R=\infty$, $\Psi^{\infty}=\Phi^{\infty}=\partial_t$
 which is consistent  with $\lim_{R\rightarrow \infty}\omega^R=0.$ Hence by definition
   $E^{\infty}=0.$

 To see the   continuity of $E$ at $R=\infty$,
 note that

 $$\|E_+(u, R)\|_{k-1, p, \delta}\leq \|\tau_{-R} ( \beta_-\beta_+'/D)\tau_{-2R} {u}_{-}\|_{k-1, p, \delta} +
 \|\tau_{-R}(\beta_- \beta_-'/D)u_+\|_{k-1, p, \delta}$$ $$ \leq C\{\|\tau_{-R}\beta_+'\tau_{-2R} {u}_{-}\|_{k-1, p, \delta}+\|\beta_-'\|_{C^{k-1}}\cdot \|\tau_{-R}(\beta_-/D)u_+\|_{k-1, p, \delta}\}$$ $$ \sim e^{2\delta (-d+l)}\|\beta_+'\|_{C^{k-1}}\cdot \| {u}_{-}\|_{k-1, p,\delta}+1/L(R)\| {u}_{+}\|_{k-1, p,\delta},$$ which goes to zero as $R\rightarrow \infty$ as long as $\|u\|_{k, p,\delta}\leq M$  for some fixed $M>>0.$ Similarly for $\|E_-(u, R)\|_{k-1, p, \delta}$ so that
 $\|E(u, R)\|_{k-1, p, \delta}\leq 1/L(R)\| {u}\|_{k, p,\delta}.$

 Consequently,
 $$\|E(u, R)-E(u_0, \infty)\|_{k-1, p, \delta}=\|E(u, R)\|_{k-1, p, \delta}\leq 1/L(R)\| {u}\|_{k, p,\delta}$$ $$ \leq 1/L(R)\{\| {u}-u_0\|_{k, p,\delta}+\| u_0\|_{k, p,\delta}\}.$$ This implies the continuity of $E$ at $R=\infty.$

 Here we  have  used the estimate
 $$\|\tau_{-R}\beta_+'\tau_{-2R} {u}_{-}\|_{k-1, p, \delta}$$
 $$\leq e^{2\delta (-d+l)}\|\beta_+'\|_{C^{k-1}}\cdot \| {u}_{-}\|_{k-1, p,\delta}$$  proved below in this section.

 \QED

  %All the results above for $\Psi^R_L$ are true for  $\Psi^a_L$  with the same proofs by replacing $``R''$ by $``a''$.

  Let ${ W}={ W}_{k, \delta}^p$  be  a small neighborhood of the space of $L_{k,\delta}^p
  $-maps with the domain $S$ near the initial ma p $f$.
  To compute the derivative

  \noindent
  $ (D_{ W}\Psi^R_L)_{(u_-, u_+)}((\xi_-, \xi_+)$ at $(u_-, u_+)$ of  $\Psi^R_L$ at $u=(u_-, u_+)$ in the direction $\xi=(\xi_-, \xi_+)$ for $\xi\in T_{(u_-, u_+)}{ W}$, we identify  the tangent space $T_{(u_-, u_+)}W=L_{k, \delta}^p(S, u^*TM)$ at $u=(u_-, u_+)$ with $T_{(f_-, f_+)}W=L_{k, \delta}^p(S, f^*TM)$ of the initial map $f$ by the usual trivialization of $T{ W}$.
  In fact in our case,   ${ W}$ is a small neighborhood of the space of $L_{k,\delta}^p
$-maps with the domain $S$ near $f$. Since the images of all such  maps lying a small neighborhood $U(f(d))$ the  double point $f(d)$, which can be identified with a small ball $B\subset E={\bf C}^n$, we may consider ${ W}$ as an open subset  of   $L_{k, \delta}^p(S, E)$ so that $T{ W}\simeq { W}\times L_{k, \delta}^p(S, E)$. Thus  $\xi=(\xi_-, \xi_+)$ can be thought as an element in $L_{k, \delta}^p(S, E)$ independent of $u$.

  Then we have
  \begin{lemma}
  $$ (D_{ W}\Psi^{R}_L)_{(u_-, u_+)}(\xi_-, \xi_+)=\partial_t (\xi_-, \xi_+)+(D_{ W}E)_{(u_-, u_+)}((\xi_-, \xi_+)$$ $$  =\partial_t (\xi_-, \xi_+)+$$ $$(
  \tau_{R}( \beta_+\beta_+'/D) \xi_-   + \tau_{R} (\beta_+\beta_-'/D){\tau_{2R}	\xi_{+}},
  \tau_{-R} ( \beta_-\beta_+'/D)\tau_{-2R} {\xi}_{-} + \tau_{-R}(\beta_- \beta_-'/D)\xi_+).$$

   $$ (\partial _{R}\Psi^R_L)_{(u_-, u_+)}=(\partial _{R}E)_{(u_-, u_+)}=$$

   $$ (\{\partial _{R}\tau_{R}( \beta_+\beta_+'/D)\} u_-   + \{\partial _{R}\tau_{R} (\beta_+\beta_-'/D)\}{\tau_{2R}	 u_{+}}+\tau_{R} (\beta_+\beta_-'/D)\partial _{R}{\tau_{2R}	u_{+}},$$ $$
\{\partial _{R}\tau_{-R} ( \beta_-\beta_+'/D)\}\tau_{-2R} {u}_{-}+ \tau_{-R} ( \beta_-\beta_+'/D)\partial _{R}\tau_{-2R} {u}_{-}+\{ \partial _{R}\tau_{-R}(\beta_- \beta_-'/D)\}u_+).$$

 $$(\partial _{\theta}\Psi^{R_\theta}_L)_{(u_-, u_+)}
 =(\partial_{\theta}E)_{(u_-, u_+)}=$$
$$ (\tau_{R_{\theta}} (\beta_+\beta_-'/D)\partial_{\theta}\tau_{2R_{\theta}}
		u_{+},
  \tau_{-R_{\theta}} ( \beta_-\beta_+'/D)\partial _{{\theta}}\tau_{-2R_{\theta}} {u}_{-}).$$
  \end{lemma}

  Here $\pm R_{\theta} =(\pm R, \pm \theta) $ and the action of $\tau_{\pm R_{\theta}}$ on $[0, \pm \infty)\times S^1$   is $(t, s)\rightarrow (t\pm R, s\pm\theta).$ We will use this kind of notations in the rest of the paper.

   The main theorem is

   \begin{theorem}
  The section  	$\Psi_L$ is of class $C^1$.
  \end{theorem}

   The proof of the theorem will be divided into two parts according to if $R=\infty.$
 
 In order to prove next lemma we list a few general facts  that will be used repeatedly:

 \medskip
 \noindent
 (A) $F_+:{W}\times [0, \infty)
 \rightarrow L(L_{k,\delta}^p(C_+, E), L_{k-1,\delta}^p(C_+, E))$
 defined by
 $F(u, R)(\xi)=\tau_{ R}\xi$ is continuous. There is a corresponding statement for the function 
 $F_-.$
 
 \medskip
 \noindent
 (B) Any smooth function $f_{\pm }$ such as $f_{\pm }=\beta_{\pm}':C_{\pm }\rightarrow {\bf R}^1$ gives rise
 a $C^{\infty}$-map $F_{ \pm}: {\bf R}^1 \rightarrow  C^m (C_{\pm},{\bf R}^1)$ defined by  $F_{\pm}(R)=f\circ \tau_R$ for any $m$. In particular, we may assume that  $m>>k$.

 \medskip
 \noindent
 (C) The paring
 
 \noindent
 $L_{k}^p(C_{\pm}, E)\times L(  L_{k, \delta }^p(C_{\pm}, E), L_{k, \delta }^p(C_{\pm}, E))\rightarrow L_{k, \delta }^p(C_{\pm}, E)$ is bounded bilinear and hence smooth as long as the space $L_{k }^p(C_{\pm}, E))$ forms  Banach algebra.
 
 \medskip
 \noindent
 (D)
 \noindent
For $m>>k$, $L(  L_{k, \delta }^p(C_{\pm}, E), L_{k, \delta }^p(C_{\pm}, E))$ is a $C^m(C_{\pm}, {\bf R}^1)$-module and  the multiplication map    $$C^m(C_{\pm}, {\bf R}^1)\times L(  L_{k, \delta }^p(C_{\pm}, E), L_{k, \delta }^p(C_{\pm}, E))\rightarrow L(  L_{k, \delta }^p(C_{\pm}, E), L_{k, \delta }^p(C_{\pm}, E))$$ is bounded bilinear and hence smooth.
 
  The proofs for  (B) and (D) are straightforward and (C)
  is stated  in the first chapter of Lang's book \cite{4}.
 The property (A)  was proved in the  last section of this paper. 
 
   \begin{lemma}
   For $R\not = \infty$, $ \Psi_L$ is of class $C^1$. Moreover for  $\Psi_L=M+E$  with $M=\partial_t, $  $D_WM$ is  continuous  as a map $D_{ W}M: { W}\times (R_0, \infty]\rightarrow
   L(L_{k, \delta}^p(S, E), L_{k-1, \delta}^p(S, E)).$  Here $W\subset L_{k, \delta}^p(S, E)$ is a small neighborhood of the initial map $f$ in  $ L_{k, \delta}^p(S, E).$
   \end{lemma}

   \proof

   The second statement is clear since $M=\partial_t  $ is linear and bounded so that $D_{ W}M_{(u_-, u_+, R)}  (\xi_-,\xi_+)=\partial _t(\xi_-,\xi_+).$
Then
   considered as a map,  $D_{ W}M: { W}\times (R_0, \infty]\rightarrow
   L(L_{k, \delta}^p(S, E), L_{k-1, \delta}^p(S, E))$ is a constant map, hence
    continuous. Here $L(L_{k, \delta}^p(S, E), L_{k-1, \delta}^p(S, E))$  is the space of bounded, hence continuous linear maps between the Banach spaces 
   $L_{k, \delta}^p(S, E)$ and $ L_{k-1, \delta}^p(S, E)
   $ under the operator norm.

 By the formulas for $D\Psi_L$ above and the  properties listed  above,  the first statement can be proved by a modified version of (A) stated in Lemma 4.4. 
  Note that  the second component of the error term of  
 $ (D_{ W}\Psi^{R}_L)_{(u_-, u_+)}(\xi_-, \xi_+)$ is $
 \tau_{-R} ( \beta_-\beta_+'/D)\tau_{-2R} {\xi}_{-} + \tau_{-R}(\beta_- \beta_-'/D)\xi_+$ with the domain with  $t$-range $[0, \infty)$ (in the natural coordinate $(t_+, s_+)$). We only need to deal with the first term
  $\tau_{-R} ( \beta_-\beta_+'/D)\tau_{-2R} {\xi}_{-}$. Then the $t$-range of $\xi_-$ is $(-\infty, 0]$ so that $\tau_{-2R} {\xi}_{-}=\xi\circ \tau_{-2R}$, which is supposed to have positive $t$-range of the domain, can only defined on $[0, 2R]$. However, since  the support of $\ \beta_+'$ is $[-d-l, -d+l]$   the support of $\tau_{-R} \beta_+'= \beta_+'\circ \tau_{-R}$ is $[R-d-l, R-d+l]$  so that  $\tau_{-R} ( \beta_-\beta_+'/D)\tau_{-2R} {\xi}_{-}$ becomes a well-defined function on $[0, \infty)\times S^1.$   Thus we  consider the function 
   $F: [0, \infty)
  \rightarrow L(L_{k,\delta}^p(C_-, E), L_{k-1,\delta}^p(C_+, E))$
  defined by
  $F( R)(\xi)=\tau_{-R} \beta_+'\cdot \tau_{ -2R}\xi_-$.
  This is a special case considered in Lemma 4.4 with $f_R$ there  being $\beta'_+$ so that   $F$  
  is continuous by Lemma 4.4. 
  
  %Another way of doing this  is to show that 
 % $$G:[0, \infty)
 % \rightarrow L(L_{k,\delta}^p(C_-, E), L_{k-1,\delta}^p((0,2R)\times S^1, E))$$
%  defined by
 % $F( R)(\xi)= \tau_{ -2R}\xi_-$ is continuous by essentially the  same argument of the proof for (A).

 \QED

 Next  we need to show that $DE$ and hence $D\Psi_L$ can be extended continuously over $R=\infty.$ This will be established by the estimates in the rest of this section.

  For $E^R(u_-, u_+)=(E^R_-(u_-, u_+), E_+^R(u_-, u_+)), $ the two components are of the same natural.
 We will only  consider $E_+^R=(\tau_{-R} ( \beta_-\beta_+'/D)\tau_{-2R} {u}_{-}-\tau_{-R}(\beta_- \beta_-'/D)u_+ )$.

 The second term
 $\tau_{-R}(\beta_- \beta_-'/D)u_+$ without involving the actions of the translation operators behaves    as expected.  The main technique reason for intruding the connection matrix $\omega$ is to get a well-behaved term
 $\tau_{-R} ( \beta_-\beta_+'/D)\tau_{-2R} {u}_{-}$
  here.

 %\begin{lemma}
 %On the interval $(R-d-l,R-d+l)$ of length $2l$ where $\tau_{-R}\beta_+'\not=0$, the weight function $e(t)$ of $ \tau_{-2R} {u}_{-}$  satisfies the bounds $e^{\delta (R-d-l)}\leq e(t)\leq e^{\delta (R-d+l)}$ so that
 %$$\|\tau_{-R}\beta_+'\tau_{-2R} {u}_{-}\|_{0, p, \delta}$$
 %$$\leq e^{2\delta (-d+l)}\|\beta_+'\|_{C_0}\cdot \| {u}_{-}\|_{0, p,\delta}.$$

 %Similarly, $$\|\tau_{-R}\beta_+'\partial_{R}(\tau_{-2R} {u}_{-})\|_{0, p, \delta}=\|\tau_{-R}\beta_+'\tau_{-2R} {u}'_{-}\|_{0, p, \delta}$$
 %$$\leq e^{2\delta (-d+l)}\|\beta_+'\|_{C_0}\cdot \| {u}'_{-}\|_{0, p,\delta}= e^{2\delta (-d+l)}\|\beta_+'\|_{C_0}\cdot \| {u}_{-}\|_{1, p,\delta}.$$
 %\end{lemma}

 %\proof

 % $$\|\tau_{-R}\beta_+'\tau_{-2R} {u}_{-}\|_{0, p, \delta}=\|(\tau_{-R}\beta_+'\tau_{-2R} {u}_{-})|_{(R-d-l,R-d+l)}\|_{0, p, \delta}$$ $$\leq e^{\delta (R-d+l)}\|\beta_+'\|_{C_0}\cdot \| {u}_{-}(t-2R))|_{t\in(R-d-l,R-d+l)}\|_{0, p}$$
 %	$$\leq e^{\delta (R-d+l)}\|\beta_+'\|_{C_0}\cdot \| {u}_{-}(v)|_{v\in(-R-d-l,-R-d+l)}\|_{0, p}$$
 %	$$\leq e^{\delta (R-d+l)}\cdot e^{(\delta (-R-d+l))}\|\beta_+'\|_{C_0}\cdot \| e(v){u}_{-}(v)|_{v\in(-R-d-l,-R-d+l)}\|_{0, p}.$$
 	%$$\leq e^{\delta (R-d+l)}\cdot e^{(\delta (-R-d+l))}\|\beta_+'\|_{C_0}\cdot \| {u}_{-}\|_{0, p,\delta}.$$

  % \QED

   \begin{lemma}
   	On the interval $(R-d-l,R-d+l)$ of length $2l$ where $\tau_{-R}\beta_+'\not=0$, the weight function $e(t)$ of $ \tau_{-2R} {u}_{-}$  satisfies the bounds $e^{\delta (R-d-l)}\leq e(t)\leq e^{\delta (R-d+l)}$ so that
   	$$\|\tau_{-R}\beta_+'\tau_{-2R} {u}_{-}\|_{k-1, p, \delta}$$
   	$$\leq e^{2\delta (-d+l)}\|\beta_+'\|_{C^{k-1}}\cdot \| {u}_{-}\|_{k-1, p,\delta}.$$
   	
   	Similarly, $$\|\tau_{-R}\beta_+'\partial_{R}(\tau_{-2R} {u}_{-})\|_{k-1, p, \delta}=\|\tau_{-R}\beta_+'\tau_{-2R} {u}'_{-}\|_{k-1, p, \delta}$$
   	$$\leq  e^{2\delta (-d+l)}\|\beta_+'\|_{C^{k-1}}\cdot \| {u}_{-}\|_{k, p,\delta}.$$
   	
   	and
   	$$\|\tau_{-R}\beta_+'(\tau_{-2R} {\xi}_{-})\|_{k-1, p, \delta}$$
   	$$\leq e^{2\delta (-d+l)}\|\beta_+'\|_{C^{k-1}}\cdot \| \xi_{-}\|_{k-1, p,\delta}.$$
   	
   \end{lemma}

\proof
    
    $$\|\tau_{-R}\beta_+'\tau_{-2R} {u}_{-}\|_{k-1, p, \delta}=\|(\tau_{-R}\beta_+'\tau_{-2R} {u}_{-})|_{(R-d-l,R-d+l)}\|_{k-1, p, \delta}$$ $$\leq e^{\delta (R-d+l)}\|\beta_+'\|_{C^{k-1}}\cdot \| {u}_{-}(t-2R))|_{t\in(R-d-l,R-d+l)}\|_{k-1, p}$$
    $$\leq e^{\delta (R-d+l)}\|\beta_+'\|_{C^{k-1}}\cdot \| {u}_{-}(v)|_{v\in(-R-d-l,-R-d+l)}\|_{k-1, p}$$
    $$\leq e^{\delta (R-d+l)}\cdot e^{(\delta (-R-d+l))}\|\beta_+'\|_{C^{k-1}}\cdot \| e(v){u}_{-}(v)|_{v\in(-R-d-l,-R-d+l)}\|_{k-1, p}.$$
    $$\leq e^{\delta (R-d+l)}\cdot e^{(\delta (-R-d+l))}\|\beta_+'\|_{C^{k-1}}\cdot \| {u}_{-}\|_{k-1, p,\delta}.$$

  \QED

   \begin{lemma}

   The interval	where $\partial _{R}\tau_{-R} ( \beta_-\beta_+'/D)\}\tau_{-2R} {u}_{-}\not=0$ is the same  $(R-d-l,R-d+l)$ so that
   	$$\|\{\partial _{R}\tau_{-R} ( \beta_-\beta_+'/D)\}\tau_{-2R} {u}_{-}\|_{k-1, p, \delta}$$ $$\leq e^{2\delta (-d+l)}\|\beta_+'\|_{C^{k-1}}\cdot \| {u}_{-}\|_{k-1, p,\delta}.$$
   \end{lemma}

   \proof 	

   By the proof of  lemma before, 	 we only need to prove the first statement.

   $$\partial _{R}\tau_{-R}(\beta_- \beta_+'/D)=
  \partial _{R}\{(\beta_- \beta_+'/D)\circ \tau_{-R}\}= \partial _{R}\{(\beta_- \beta_+'/D)(t-R)\}$$ $$=-((\beta_- \beta_+'/D)'\circ \tau_{-R} =-\{(\beta_- /D)'\beta_+'+(\beta_- /D)\beta_+''\}\circ \tau_{-R}$$ $$ =\tau_{-R}(\beta_- /D)'\cdot \tau_{-R}\beta_+'+\tau_{-R} (\beta_- /D)\cdot\tau_{-R}\beta_+''.$$

  Since outside $(R-d-l,R-d+l)$, $\tau_{-R}\beta_+'=\tau_{-R}\beta_+''=0,$  so is $\partial _{R}\tau_{-R}(\beta_- \beta_+'/D).$

 \QED

   Hence

   \begin{lemma}
   	$ \|(D_{ W}\Psi_L)_+(u_-, u_+, R)-\partial_t\|_{o}\rightarrow 0$ as $R\rightarrow \infty$ uniformly in $u$.
\end{lemma}

   \proof

   	$$ \|(D_{ W}\Psi_L)_+(u_-, u_+, R)-\partial_t\|_{o}=\sup_{\|\xi||_{k, p, \delta}\leq 1}\|(D_{ W}\Psi_L)_+(u_-, u_+, R)(\xi)-\partial_t(\xi_+)\|_{k-1, p, \delta}$$
   	$$=\sup_{\|\xi||_{k, p, \delta}\leq 1}\|(D_{ W}\Psi^{R}_{L, +})_{(u_-, u_+)}(\xi_-, \xi_+)-\partial_t(\xi_+)\|_{k-1, p, \delta} $$ $$\leq \sup_{\|\xi||_{k, p, \delta}\leq 1}\{\|\tau_{-R} ( \beta_-\beta_+'/D)\tau_{-2R} {\xi}_{-}\|_{k-1, p, \delta} + \|\tau_{-R}(\beta_- \beta_-'/D)\xi_+\|_{k-1, p, \delta} \}$$

   	$$\leq \sup_{\|\xi||_{k, p, \delta}\leq 1}\{||\beta_-/D||_{C^{k-1}}\|\tau_{-R} ( \beta_+')\tau_{-2R} {\xi}_{-}\|_{k-1, p, \delta} + \|\beta_- \beta_-'/D\|_{C^{k-1}}\|\xi_+\|_{k-1, p, \delta} \}$$
   	
   	$$\leq \sup_{\|\xi||_{k, p, \delta}\leq 1}\{e^{2\delta (-d+l)}||\beta_-/D||_{C^{k-1}}\|\beta_+'\|_{C^{k-1}}\cdot \| \xi_{-}\|_{k-1, p,\delta} + \|\beta_- \beta_-'/D\|_{C^{k-1}}\|\xi_+\|_{k-1, p, \delta} \}.$$
   	
   Now the key point is $$\|\beta_- \beta_-'/D\|_{C^{k-1}}\sim
   \|\beta_- \beta_-'/D\|_{C^{0}}\sim \|\beta_- /D\|_{C^{0}}\|\beta_-'\|_{C^{0}}\sim 1/L(R)\rightarrow 0$$ as $R\rightarrow \infty.$

   The  first $``\sim''$ above  can be proved inductively. Indeed
   since $$(\beta_- \beta_-'/D)'=\{((\beta_-')^2+\beta_- \beta_-'')D-(\beta_- \beta_-')D'\}/D^2$$ $$ =\{((\beta_-')^2+\beta_- \beta_-'')D-2(\beta_- \beta_-')(\beta_- \beta_-'+\beta_+ \beta_+')\}/D^2,$$

   it is easy to see  that  for any $i>0$, each term of $\nabla^i\{(\beta_- \beta_-'/D)\}$  is a product of  the  terms of the form containing at least one of
     $\beta_-' $ and  $\beta_+'$,   or   their  derivatives. Hence  $||\{(\beta_- \beta_-'/D)\}||_{C^0}$ is  the  lowest order term in $1/L(R)$  inside
   $ \|(\beta_- \beta_-'/D)\|_{C^{k-1}}.$

    This proves that $$ \|(D_{ W}\Psi_L)_+(u_-, u_+, R)-\partial_t\|_{o}
   \leq ||\beta_-/D||_{C^{k-1}}\{e^{2\delta (-d+l)}\|\beta_+'\|_{C^{k-1}}+\|\beta_-'||_{C^{0}}\} \sim 1/L(R)\rightarrow 0$$ as $R\rightarrow \infty.$

   \QED

   \begin{lemma}
   	$$ \|(\partial _{r}(\Psi^R_L)_+)_{(u_-, u_+)}\|_{k-1, p,\delta}\sim 1/(\ln^2 R)\{\| {u}_{-}\|_{k, p,\delta}+\| {u}_{+}\|_{k, p,\delta}\}$$

   \end{lemma}

   	\proof
   	
   	$$ \|(\partial _{R}(\Psi^R_L)_+)_{(u_-, u_+)}\|_{k-1, p,\delta}
   	=\|(\partial _{R}E_+)_{(u_-, u_+)}\|_{k-1, p, \delta}\leq $$
   	
   	  $$
   	\|\{\partial _{R}\tau_{-R} ( \beta_-\beta_+'/D)\}\tau_{-2R} {u}_{-}\|_{k-1, p, \delta}+ \|\tau_{-R} ( \beta_-\beta_+'/D)\partial _{R}\tau_{-2R} {u}_{-}\|_{k-1, p, \delta}+$$ $$ \|\{ \partial _{R}\tau_{-R}(\beta_- \beta_-'/D\}u_+)\|_{k-1, p, \delta} $$

   $$\leq e^{2\delta (-d+l)}\cdot \|( \beta_-\beta_+'/D)||_{C^k}\| {u}_{-}\|_{k-1, p,\delta}
   +e^{2\delta (-d+l)}
   \cdot \| \beta_-/D||_{C^{k-1}}\|\beta_+'\|_{C^{k-1}}\cdot
   	\| {u}_{-}\|_{k, p,\delta}$$ $$ + \|\{ \partial _{R}\tau_{-R}(\beta_- \beta_-'/D)\}\|_{C^{k-1}}\cdot \|u_+\|_{k-1, p, \delta}$$

    $$\leq  e^{2\delta (-d+l)}\cdot C(\|\beta\|_{C^{k+1}})\| {u}_{-}\|_{k, p,\delta}+ \|(\beta_- \beta_-'/D)'\circ\tau_{-R}\|_{C^{k-1}}\cdot \|u_+\|_{k-1, p, \delta}$$

 Now  $ \|(\beta_- \beta_-'/D)'\|_{C^{k-1}}\sim
 \|(\beta_- \beta_-'/D)'\|_{C^{0}}$ and  the key point is $$\|(\beta_- \beta_-'/D)'\circ\tau_{-R}\|_{C^{k-1}}=\|(\beta_- \beta_-'/D)'\|_{C^{k-1}}$$ $$ \sim\|(\beta_- \beta_-'/D)'\|_{C^{0}}\sim  \|\{((\beta_-')^2+\beta_- \beta_-'')D-(\beta_- \beta_-')D'\}/D^2)\|_{C^0}$$ $$\sim ||\{((\beta_-')^2+\beta_- \beta_-'')D-(\beta_- \beta_-')(\beta_-\beta_-'+\beta_+\beta_+')\}/D^2)\|_{C^0}$$
 $$\sim ||(\beta_-')^2\|_{C^0}+|| \beta_-''\|_{C^0}+|| \beta_-'\beta_+'\|_{C^0}\sim 1/L(R)^2=1/(R\{\ln R\}^4);$$
  
   This first  statement can be proved inductively. Indeed  it is easy to see  from above  that  for any $i>0$, each term of $\nabla^i\{(\beta_- \beta_-'/D)'\}$  is a product of  the  terms containing at least one of the  terms
   $(\beta_-')^2$, $\beta_-''$,  $\beta_-'\cdot\beta_+'$  or   their  derivatives so that  for  $i>0$, $||\nabla^i\{(\beta_- \beta_-'/D)'\}||_{C^0}$ is a lower order term comparing with
   $ \|(\beta_- \beta_-'/D)'\|_{C^{0}}.$

  Thus 	$$ \|(\partial _{R}(\Psi^R_L)_+)_{(u_-, u_+)}\|_{k-1, p,\delta}\sim  1/(R\{\ln R\}^4)\{\| {u}_{-}\|_{k, p,\delta}+\| {u}_{+}\|_{k, p,\delta}\}$$ so that
  $$ \|(\partial _{r}(\Psi^R_L)_+)_{(u_-, u_+)}\|_{k-1, p,\delta}=\|(\partial _{R}(\Psi^R_L)_+)_{(u_-, u_+)}\|_{k-1, p,\delta}|dR/dr|
 $$ $$ \sim  1/(R\{\ln R\}^4)\{\| {u}_{-}\|_{k, p,\delta}+\| {u}_{+}\|_{k, p,\delta}\}\cdot |dR/dr|$$ $$= 1/(R\{\ln R\}^4)\{\| {u}_{-}\|_{k, p,\delta}+\| {u}_{+}\|_{k, p,\delta}\}R\cdot (\ln R)^2=1/(\ln^2 R)\{\| {u}_{-}\|_{k, p,\delta}+\| {u}_{+}\|_{k, p,\delta}\}$$

 since $R=e^{1/r}$ and ${\frac {dR}{dr}}={\frac {de^{1/r}}{dr}=e^{1/r}\cdot (-1/r^2) }=-R\cdot (\ln R)^2.$

\QED

  \begin{lemma}

  		$$(\partial _{\theta}(\Psi^{R_\theta}_L)_+)_{(u_-, u_+)}\sim Ce^{\delta (-2d+2l)}\cdot  \|{u}_{-} \|_{k, p, \delta}.$$
  	Hence 	$$1/r(\partial _{\theta}(\Psi^{R_\theta}_L)_+)_{(u_-, u_+)}<<e^{-1/r} \|{u}_{-} \|_{k, p, \delta}$$
  \end {lemma}

 \proof

 $$(\partial _{\theta}(\Psi^{R_\theta}_L)_+)_{(u_-, u_+)}
 =\|(\partial_{\theta}E)_{(u_-, u_+)}\|_{k-1, p,\delta}=$$

 $$
 \|	\tau_{-R_{\theta}} ( \beta_-\beta_+'/D)\partial _{{\theta}}\tau_{-2R_{\theta}} {u}_{-}\|_{k-1, p,\delta}.$$

 Recall  that $
 	\tau_{-R_{\theta}} ( \beta_-\beta_+'/D)\not=0$ only when $t\in (R-d-l,R-d+l)$ and on this interval the weight function $e(t)$ of $ \tau_{-2R} {u}_{-}$  satisfies the bounds $e^{\delta (R-d-l)}\leq e(t)\leq e^{\delta (R-d+l)}$.

 Thus on this interval
 $$
 \|	\tau_{-R_{\theta}} ( \beta_-\beta_+'/D)\partial _{{\theta}}\tau_{-2R_{\theta}} {u}_{-}\|_{k-1, p,\delta}\sim C\|\partial _{{\theta}}\tau_{-2R_{\theta}} {u}_{-}|_{[R-d-l,R-d+l]}\|_{k-1, p,\delta}$$
 $$ \sim Ce^{\delta (R-d+l)}\cdot \|\partial _{s}{u}_{-}
 (t-2R, s-2\theta)|_{[R-d-l,R-d+l]}\|_{k-1, p}$$ $$ =
 Ce^{\delta (R-d+l)}\cdot \|\partial_{s}{u}_{-}(t, s) |_{[-R-d-l,-R-d+l]}\|_{k-1, p}$$
 $$=Ce^{\delta (R-d+l)}\cdot e^{\delta(-R-d+l)}
  \|e_-(t)\partial_{s}{u}_{-}(t, s) |_{[-R-d-l,-R-d+l]}\|_{k-1, p}$$
   $$=Ce^{\delta (-2d+2l)}\cdot  \|\partial_{s}{u}_{-} \|_{k-1, p, \delta}=Ce^{\delta (-2d+2l)}\cdot  \|{u}_{-} \|_{k, p, \delta}.$$

 	Hence 	$$1/r(\partial _{\theta}(\Psi^{R_\theta}_L)_+)_{(u_-, u_+)}\sim C\ln R \cdot  e^{\delta (-d)}\cdot  \|{u}_{-} \|_{k, p, \delta}.$$ Recall $d\sim l\sim R^{1/2}$ so that $\ln R \cdot e^{\delta (-d)}\|{u}_{-} \|_{k, p, \delta}\sim e^{-\delta R^{1/2}}\|{u}_{-} \|_{k, p, \delta}<< 1/R\|{u}_{-} \|_{k, p, \delta}=e^{-1/r}\|{u}_{-} \|_{k, p, \delta}.$

 \QED

 \begin{pro}
 Away from $r=0,$ the partial derivatives $\partial _x \Psi_L$ and $\partial _y \Psi_L$ exits and continuous, and they can be extended continuously over
 $r=0.$
 \end{pro}

\proof

This follows from the estimates established so far together with  the following elementary formula:

$$
\left[\begin{array}{c}
{\frac {\partial}{\partial x}}\\
{\frac {\partial}{\partial y}}
\end{array}\right]
=\left[\begin{array}{ll}
\cos \theta  & -\sin \theta \\
\sin \theta & \cos \theta
\end{array}\right]
\left[\begin{array}{c}
{\frac {\partial}{\partial r}}\\
	1/r\cdot {\frac {\partial}{\partial\theta}}
\end{array}\right].
$$

Indeed, by this formula
 $\partial_x\Psi_L:L_{k, \delta}^p(C_-, E)\times  L_{k, \delta}^p(C_+, E)\times D^*_{r_0}\rightarrow L_{k-1, \delta}^p(C_-, E)\times  L_{k-1, \delta}^p(C_+, E)$ is given by $$\partial_x\Psi_L(u, (x, y))=\partial_x\Psi_L(u, (R, \theta))$$ $$ =\cos \theta \partial _r   
\Psi_L (u, (R, \theta))-\sin\theta \cdot 1/r\partial _{\theta}\Psi_L (u, (R, \theta)).$$

The proof for the case $k=1$ is clear. For the general $k$, note that
For $i+j\leq k-1,$
$$D_s^iD_t^j \{\partial_x\Psi_L(u, (R, \theta))\}$$ $$ =\cos \theta \{D_s^iD_t^j \partial _r   
\Psi_L (u, (R, \theta))\}-\sin\theta \cdot \{D_s^iD_t^j 1/r\partial _{\theta}\Psi_L (u, (R, \theta))\}$$ so that $$\| \partial_x\Psi_L(u, (R, \theta))\|_{k-1, p, \delta}$$ $$ \leq \|\partial _r   
\Psi_L (u, (R, \theta))\|_{k-1, p, \delta}+|| 1/r\partial _{\theta}\Psi_L (u, (R, \theta))\|_{k-1, p, \delta}.$$

Then the conclusion follows.

%We leave the general case
%using $L_{k, \delta}^p$-norm to the readers since we do not need it.

%Indeed it is clear without any further argument  that  for each fixed $\theta_0$,  $\partial _x \Psi_L$ and $\partial _y \Psi_L$ extends over $r=0$ along the line $\theta=\theta_0$ even using the general $L_{k, \delta}^p$-norm. In particular $\partial _x \Psi_L$  extends over $x=0$ along $x$-axis. Similarly 
% $\partial _y \Psi_L$  extends over $y=0$ along $y$-axis.  By the last lemma of the next section, this  is already  sufficient  for proving that $\Psi_L$ is of class $C^1$.

\QED
\section{Some basic estimates}

   In this last section,  we collect  some basic estimates  used in paper first, then use them to finish the proof of the main theorem of this paper.
   The proofs are only given in $L^p$-norm with $p>1$. The corresponding results for $L^p_k$-norm can be derived from this by  replacing a $L_k^p$-function $f$ by a $L^p$-function $(f, Df, \cdots, D^kf).$ For the purpose of this paper,
   we assume in addition that $1-2/p>0$.

   The following inequality will be used repeatedly.

   \begin{lemma}
   Let $ F(x, t)$ be a smooth function for $(x, t)\in \Sigma\times [0, 1]$ with compact support. Then $\int_{\Sigma}|\int_{[0, 1]} F(x, t)dt|^pd vol_{\Sigma}\leq  \int_{[0, 1]}\int_{\Sigma} | F(x, t)|^pd vol_{\Sigma}dt.$
   \end{lemma}

 \proof

  Since the functions $f_1(x)=|x|^p$  is convex, for each fixed $x$ and $a$ and $b$,
  $ |t F(x, a)+(1-t)F(x, b)|^p\leq t| F(x, a)|^p+(1-t)|F(x, b)|^p$ so that
   $|\Sigma_{i=1}^nF(x, i/n)/n|^p\leq \Sigma_{i=1}^n|F(x, i/n)|^p/n.$

  Hence  $$|\int_{[0, 1]} F(x, t)dt|^p\leq \int_{[0, 1]} |F(x, t)|^pdt$$  and
   $$\int_{\Sigma}|\int_{[0, 1]} F(x, t)dt|^pdvol_{\Sigma}\leq \int_{\Sigma}\int_{[0, 1]} |F(x, t)|^pdtdvol_{\Sigma}\leq  \int_{[0, 1]}\int_{\Sigma} | F(x, t)|^pd vol_{\Sigma}dt.$$

   \QED

 In our case, $\Sigma={\bf R}^1\times S^1$ which is not compact. Then the set of smooth functions with compact support is dense in the space of  $L_k^p$ or $L_{k, \delta}^p$-functions so that above is applicable.

   \begin{lemma}  	
   	$$\|\tau_{R}\xi\|_{0, p, \delta}=
   	\leq e^{-\delta R}\|\xi\|_{0, p,\delta}$$
   	for $\xi:[0, \infty)\rightarrow E$.
      	
   \end{lemma}

   \proof

   $$\|\tau_{R}\xi\|_{0, p, \delta}=\|\xi(t+R)\|_{0, p, \delta}
   =\|e(t)\xi(t+R)\|_{0, p}=e^{-\delta R}\|e(u)\xi(u)|_{u\in [R, \infty)}\|_{0, p}$$ $$ \leq e^{-\delta R}\|\xi\|_{0, p,\delta}.$$

   \QED

  % More generally we need to consider $\tau_{R}(\tau_{R}(\xi)\circ \Gamma)=\xi\circ \tau_R \Gamma\tau_R$ as a function of $R$ first.

   %Recall $\tau_R\xi=: \xi\circ \tau_R.$
   %Let  $\tau_R \Gamma\tau_R=:\phi_R$ is  a family of diffeomorphisms depending on $R$.
   Let $F: L_{1, \delta}^p(C_+, E)\times [0, \infty)\rightarrow L_{0, p,\delta}(C_+, E)$ defined by $F( \xi, R)=\xi\circ \tau_R$.  It is proved blow in this section that    $\partial_R F=\xi'\circ \tau_R$. Denote $\partial_R F$ by $G$  and $\xi'$ by $\eta$.

  % Then $$\|F_1(R_1)-F_1(R_2)\|_{0, p, \delta }=\|\xi'\circ \tau_{R_1}-\phi_{R_2}'\xi'\circ \tau_{R_2}\|_{0, p, \delta}$$ $$= \|\eta\circ \tau_{R_1}-\eta\circ \tau_{R_2}\|_{0,p, \delta}$$

   %Clearly the factor $\phi_R'$ in $F_1(R)$ does not affect the analytic property of $F_1$.

   Consider $G: L_{0, \delta}^p (C_+, E)\times[0, \infty)\rightarrow L_{0, \delta}^p(C_+, E)$,  $G( \eta, R)=\eta\circ \phi_{R}.$

    Above functions $F$ and $G$ are repeatedly used in this paper.

   \begin{lemma}

    The function $G$ is continuous.

   \end{lemma}

   \proof

    $$\|G(\eta_1, R_1)-G(\eta_2, R_2)\|\leq \|G(\eta_1, R_1)-G(\eta_1, R_2)\|+\|G(\eta_1, R_2)-G(\eta_2, R_2)\|$$
   $$ \leq \|\eta_1\circ \tau_{R_1} -\eta_1\circ \tau_{R_2}\|+\|\eta_1\circ\tau _{R_2}-\eta_2\circ\tau_{R_2}\|$$
   $$
   \leq ||\eta_1\circ \tau_{R_1} -{\hat \eta}\circ \tau_{R_1}\|+||{ \hat\eta}\circ \tau_{R_1} -{\hat \eta}\circ \tau_{R_2}||+  ||{\hat \eta}\circ \tau_{R_2} -\eta_1\circ \tau_{R_2}||+   \|(\eta_1-\eta_2)\circ\phi_{R_2}\|$$
   	$$ \leq ||(\eta_1 -{\hat \eta})\circ \tau_{R_1}\|+||{\hat \eta}\circ \tau_{R_1} -{\hat \eta}\circ \tau_{R_2}\|+||(\eta_1 -{\hat \eta})\circ \tau_{R_2}||+\|(\eta_1-\eta_2)\circ\tau_{R_2}\|$$
   $$\leq ( e^{-\delta R_1}+ e^{-\delta R_2})(||(\eta_1 -{\hat \eta})\|+||(\eta_1 -\eta_2)\|) +||{\hat \eta}\circ \tau_{R_1} -{\hat \eta}\circ \tau_{R_2}\|.$$

   The last term
    $$||{\hat \eta}\circ \tau_{R_1} -{\hat \eta}\circ \tau_{R_2}\|$$ $$ = \|{\hat \eta}\circ \tau_{R_1}-{\hat \eta}\circ\tau_{R_2} \|_{0,p, \delta}=|{R_2}-{R_1}|\cdot\|\int_{[0,1]} {\hat \eta}'((1-s)(t+{R_1})+s(t+{R_2}))ds\|_{0, p, \delta}$$
    $$\leq|{R_2}-{R_1}|\cdot \int_{[0,1]} \|{\hat \eta}'(t+(1-s){R_1}+s{R_2}))\|_{0, p, \delta} ds$$ $$ \leq |{R_2}-{R_1}|\cdot \|{\hat \eta}'\|_{0, p, \delta}\cdot \int_{[0,1]} e^{-\delta((1-s){R_1}+s{R_2})} ds$$ $$=|{R_2}-{R_1}|\cdot \{ e^{-\delta {R_2}}-e^{-\delta {R_1}}\}/\{\delta(R_2-R_1)\}\cdot  \|{\hat \eta}'\|_{0, p, \delta}=|e^{-\delta {R_2}}-e^{-\delta {R_1}}|/\delta\cdot  \|{\hat \eta}'\|_{0, p, \delta}$$ $$ =|e^{-\delta {R_2}}-e^{-\delta {R_1}}|/\delta\cdot  \|{\hat \eta}\|_{1, p, \delta}.$$

     Given any $\epsilon>0$, we may chose a smooth ${\hat  \eta}$ with compact support such that  $\|{\hat \eta}-\eta_1\|_{0, p, \delta}<\epsilon.$ Then above estimate proves the continuity of $G$.
     \QED

    For  the function
    $F(u, R) =u\circ \tau_R$ with $u\in L_{1, \delta}^p(C_+, E)$ above, it is proved below that
      $\partial_u F=\partial_1 F:L_{1, \delta}^p(C_+, E)\times {\bf R}^1\rightarrow  L(L_{1, \delta}^p(C_+,E), L_{0, \delta}^p(C_+, E))$ given by
     $(\partial_1 F)_{(u, R)}(\eta)=\eta\circ \tau_R$. It   is $u$-independent, and  hence becomes a map, denoted by  $D: {\bf R}^1\rightarrow  L(L_{1, \delta}^p(C_+, E), L_{0, \delta}^p(C_+, E)).$

%For the applications in this paper and its companion, we need the following function similar to above: $H: {\bf R}^1\rightarrow  L(L_{1, \delta}^p(C_-, E), L_{0, \delta}^p(C_+, E)),$ given by $H(R)(\xi)= \tau_{-R} ( \beta_+')\tau_{-2R} {\xi}$. Here $\xi\in L_{1, \delta}^p(C_-, E)$ with $C_-=(-\infty, 0)\times S^1.$ Hence the $t$-range of $\xi$ is $(-\infty, 0]$ so that $\tau_{-2R} {\xi}_{-}=\xi\circ \tau_{-2R}$ can only   have positive $t$-range  on $[0, 2R]$. However since   the support of $\ \beta_+'$ is $[-d-l, -d+l]$,   the support of $\tau_{-R} \beta_+'= \beta_+'\circ \tau_{-R}$ is $[R-d-l, R-d+l]$  so that  $\tau_{-R} \beta_+'\tau_{-2R} {\xi}$ becomes a well-defined function on $C_+=[0, \infty)\times S^1.$

For the applications in this paper and its companion, we need a few functions
similar to above.   We  define  the following function that is general enough  for these applications.

  Let  $H: {\bf R}^1\rightarrow  L(L_{1, \delta}^p(C_-, E), L_{0, \delta}^p(C_+, E)),$ given by $H(R)(\xi)= \tau_{-R} f_R\tau_{-2R} {\xi}$, where $f_R:{\bf R}^1\rightarrow {\bf R}^1$  is  a  $C^{\infty}$-function smoothly depending on $R\in [R_0, \infty)$ such that (i) the support of $f_R$ is in $[A(R), B(R)]=[a_1 d(R)+a_2l(R), b_1(R)+b_2l(R)]$ with $A(R)< B(R)$ and $|A(R)|, |B(R)|<< R$; (ii) For $i\geq 0, $ $\|\partial_R^i  f_R\|_{C^k}<C_k$ independent of $R$.
  Here $\xi\in L_{1, \delta}^p(C_-, E)$ with $C_-=(-\infty, 0)\times S^1.$ Hence the $t$-range of $\xi$ is $(-\infty, 0]$ so that $\tau_{-2R} {\xi}_{-}=\xi\circ \tau_{-2R}$ can only   have positive $t$-range  on $[0, 2R]$. However since   the support of $f_R$ is $[A(R), B(R)]$,   the support of $\tau_{-R} f_R= f_R\circ \tau_{-R}$ is $[R+A(R), R+B(R)]\subset [0,2R]$  so that  $\tau_{-R}f_R\tau_{-2R} {\xi}$ becomes a well-defined function on $C_+=[0, \infty)\times S^1.$

     \begin{lemma}
     	Let $D:[R_0, \infty)\rightarrow L(L_{1, p, \delta}(C_+, E), L_{0, p, \delta}(C_+,E))$ and $H:[R_0, \infty)\rightarrow L(L_{1, p, \delta}(C_-, E), L_{0, p, \delta}(C_+,E))$ be the functions  defined above.

     	Then $$\|D(R_1)-D(R_2)\|_{o}=\sup_{\|\xi|_{1, p, \delta}\leq 1}
     	\|\tau_{R_1}\xi-\tau_{R_2}\xi\|_{0, p, \delta}\leq  |e^{-\delta R_1}-e^{-\delta R_2}|/\delta.$$
     	
   and  $$\|H(R_1)-H(R_2)\|_{o}\leq  C_0\{|R_1-R_2| e^{2\delta(|A(R_1)|+|B(R_1)|)}
   +|e^{2\delta R_1}-e^{2\delta R_2}|/\delta\cdot e^{-\delta (R_1)}\} 
   .$$
     	
     	Hence $D$ and $H$ are  continuous. Here we denote the operator norms  by $\|-\|_{o}.$
     \end{lemma}

     \proof

 We only prove the lemma for $H$ since the proof for  $D$ is easier.
  
 % More generally, $\beta'_+$  can be replaced by a function $f(R)$  with support $(d(R)-l(R), d(R)+l(R))$  (for $\beta'_+$) or $(-d(R)-l(R), d(R)+l(R))$ (for application in other cases ) such that $\|f(R)\|_{C^0}<C$ independent of$R$.

     $$\|H(R_1)-H(R_2)\|_{o}=\sup_{\|\xi||_{1, p, \delta}\leq 1}
     \|\tau_{-R_1}f_{R_1}\tau_{-2R_1}\xi-\tau_{-R_2}f_{R_2}\tau_{-2R_2}\xi\|_{0, p, \delta}$$
     $$=\sup_{\|{\hat \xi}||_{1, p, \delta}\leq 1} \|\tau_{-R_1}f_{R_1}\tau_{-2R_1}{\hat \xi}-\tau_{-R_2}f_{R_2}\tau_{-2R_2}{\hat \xi}\|_{0, p, \delta}$$ 
     
      $$=\sup_{\|{\hat \xi}||_{1, p, \delta}\leq 1}
       \|[\tau_{R_1}f_{R_1}{\hat \xi}]\circ
        \tau_{-2R_1}-[\tau_{R_2}f_{R_2}{\hat \xi}]\circ \tau_{-2R_2}\|_{0, p, \delta}.$$
      
      with $ {\hat \xi}$ being  smooth  with compact support with domain containing in $C_-=(-\infty, 0 ]\times S^1$. 
      
      %However $[\tau_{R}\beta_{+}{\hat \xi}]$ has compact support in $(-R-d-l, -R-d+l)\times S^1\subset (-\infty, 0)\times S^1\subset {\bf R}^1\times S^1$ so that the argument below  makes sense.

Then $$\|[\tau_{R_1}f_{R_1}{\hat \xi}]\circ
\tau_{-2R_1}-[\tau_{R_2}f_{R_2}{\hat \xi}]\circ \tau_{-2R_2}\|_{0, p, \delta}$$ $$ \leq \|[\tau_{R_1}f_{R_1}{\hat \xi}]\circ
\tau_{-2R_1}-[\tau_{R_1}f_{R_1}{\hat \xi}]\circ
\tau_{-2R_2}\|_{0, p, \delta}+\|[\tau_{R_1}f_{R_1}{\hat \xi}]\circ \tau_{-2R_2}-[\tau_{R_2}f_{R_2}{\hat \xi}]\circ \tau_{-2R_2}\|_{0, p, \delta}.$$
Here the $L_{0,\delta}^p$-norm is taken over $[0, \infty)\times S^1$.

Assume that $R_1\leq R_2$ with $|R_2-R_1|\leq 1.$

     Then  a similar  estimate  to  the proof of the last lemma gives
     $$\|[\tau_{R_1}f_{R_1}{\hat \xi}]\circ
     \tau_{-2R_1}-[\tau_{R_1}f_{R_1}{\hat \xi}]\circ
     \tau_{-2R_2}\|_{0, p, \delta}
     $$

    $$ \leq 2|{R_2}-{R_1}|\cdot\|\int_{[0,1]} [\tau_{R_1}f_{R_1}{\hat \xi}]'((1-s)(t-2{R_1})+s(t-2{R_2}))ds\|_{0, p, \delta}$$
    $$\leq2|{R_2}-{R_1}|\cdot \int_{[0,1]} \|[\tau_{R_1}f_{R_1}{\hat \xi}]'(t-2(1-s){R_1}
    -2s{R_2})\|_{0, p, \delta} ds$$

     $$\leq \leq2|{R_2}-{R_1}|\cdot \int_{[0,1]} \|e^{\delta t}[\tau_{R_1}f_{R_1}{\hat \xi}]'(t-2(1-s){R_1}
     -2s{R_2})\|_{0, p} ds$$
    $$=2|{R_2}-{R_1}|\cdot \int_{[0,1]} \|e^{\delta[ 2(1-s){R_1}
    	+2s{R_2}]}\cdot e^{\delta u}[\tau_{R_1}f_{R_1}{\hat \xi}]'(u)\|_{0, p} ds$$
    
  (   with 
   $u= t-2(1-s){R_1}
    -2s{R_2}\in [-R_1+A(R_1), -R_1+B(R_1)]$)

    $$=2|{R_2}-{R_1}|\cdot \int_{[0,1]} e^{\delta[ 2(1-s){R_1}
    	+2s{R_2}]}ds\cdot \| e^{\delta u}[\tau_{R_1}f_{R_1}{\hat \xi}]'(u)|_{[-R_1+A(R_1), -R_1+B(R_1)]\times S^1}\|_{0, p} $$
    
    $$=|e^{2\delta R_1}-e^{2\delta R_2}|/\delta\cdot \| e^{\delta u}\cdot \tau_{R_1} [f_{R_1}\tau_{-R_1}{\hat \xi}]'(u)|_{[-R_1+A(R_1), -R_1+B(R_1)]\times S^1}\|_{0, p} $$
    
    $$=|e^{2\delta R_1}-e^{2\delta R_2}|/\delta\cdot \|  e^{\delta (v-R_1)} [f_{R_1}\tau_{-R_1}{\hat \xi}]'(v)|_{[A(R_1), B(R_1)]\times S^1}\|_{0, p} $$
    $$\leq |e^{2\delta R_1}-e^{2\delta R_2}|/\delta\cdot \|f_{R_1}\|_{C^0} \| \tau_{-R_1} e^{\delta v} \cdot \tau_{-R_1}{\hat \xi}(v)|_{[A(R_1), B(R_1)]\times S^1}\|_{1, p} $$
    
    $$\leq |e^{2\delta R_1}-e^{2\delta R_2}|/\delta\cdot \|f_{R_1}\|_{C^0}\cdot  \|  e^{\delta v} \cdot{\hat \xi}(v)|_{[-R_1+A(R_1), -R_1+B(R_1)]\times S^1}\|_{1, p} $$
    
     $$\leq |e^{2\delta R_1}-e^{2\delta R_2}|/\delta\cdot\|f_{R_1}\|_{C^0} \cdot e^{\delta (-R_1+|A(R_1)|+|B(R_1)|) }\cdot e^{\delta (-R_1+|A(R_1)|+|B(R_1)|) }
      \|  e_-(t) \cdot{\hat \xi}(t)\|_{1, p} $$
     $$= |e^{2\delta R_1}-e^{2\delta R_2}|/\delta\cdot C_0\cdot e^{-2\delta (R_1-|A(R_1)|-|B(R_1)|)} 
     \|  {\hat \xi}(t)\|_{1, p, \delta} $$
     $$\leq ||e^{2\delta R_1}-e^{2\delta R_2}|/\delta\cdot C_0\cdot e^{-\delta (R_1)} 
     \|  {\hat \xi}(t)\|_{1, p, \delta}. $$

    %$$ \leq 2|{R_2}-{R_1}|\cdot \|{\hat \eta}'\|_{0, p, \delta}\cdot \int_{[0,1]} e^{-\delta((1-s){R_1}+s{R_2})} ds$$ 

    %$$=|{R_2}-{R_1}|\cdot \{ e^{-\delta {R_2}}-e^{-\delta {R_1}}\}/\{\delta(R_2-R_1)\}\cdot  \|{\hat \eta}'\|_{0, p, \delta}=|e^{-\delta {R_2}}-e^{-\delta {R_1}}|/\delta\cdot  \|{\hat \eta}'\|_{0, p, \delta}$$ $$ =|e^{-\delta {R_2}}-e^{-\delta {R_1}}|/\delta\cdot  \|{\hat \eta}\|_{1, p, \delta}.$$
    
    and    $$\|[\tau_{R_1}f_{R_1}{\hat \xi}]\circ \tau_{-2R_2}-[\tau_{R_2}f_{R_2}{\hat \xi}]\circ \tau_{-2R_2}\|_{0, p, \delta}$$ $$ =||e_+(t)\{[\tau_{R_1}f_{R_1}{\hat \xi}]-[\tau_{R_2}f_{R_2}{\hat \xi}]\}\circ \tau_{-2R_2}\|_{0, p}.$$
      (  with 
      $u= t-2R_2
      \in [-R_1+A(R_1)-1, -R_1+B(R_1)+1]$, 
      
      noting that $|R_2-R_1|\leq 1$ )

    $$= e^{2\delta R_2}||e^{\delta u}[\tau_{R_1}f_{R_1}-\tau_{R_2}f_{R_2}]{\hat \xi}(u)|_{[-R_1+A(R_1)-1, -R_1+B(R_1)+1]\times S^1}\|_{0, p}$$
    
    $$\leq  e^{2\delta R_2}\cdot e^{\delta(-R_1+|A(R_1)|+|B(R_1)|+2)}
    ||[\tau_{R_1}f_{R_1}-\tau_{R_2}f_{R_2}]{\hat \xi}(u)|_{[-R_1+A(R_1)-1, -R_1+B(R_1)+1]\times S^1}\|_{0, p}$$
  
  $$\leq  e^{2\delta R_2}\cdot e^{\delta(-R_1+|A(R_1)|+|B(R_1)|+2)} \|\tau_{R_1}f_{R_1}-\tau_{R_2}f_{R_2}\|_{C_0}
  ||{\hat \xi}(u)|_{[-R_1+A(R_1)-1, -R_1+B(R_1)+1]\times S^1}\|_{0, p}$$

  $$\leq  e^{2\delta R_2}\cdot e^{\delta(-R_1+|A(R_1)|+|B(R_1)|+2)} \|\tau_{R_1}f_{R_1}-\tau_{R_2}f_{R_2}\|_{C_0}\cdot e^{\delta(-R_1+|A(R_1)|+|B(R_1)|+2)}
  ||e_-{\hat \xi}(u)|_{[-R_1+A(R_1)-1, -R_1+B(R_1)+1]\times S^1}\|_{0, p}$$
  
  $$\sim    e^{2\delta(|A(R_1)|+|B(R_1)|)} \|\tau_{R_1}f_{R_1}-\tau_{R_2}f_{R_2}\|_{C_0}\cdot 
  ||{\hat \xi}\|_{0, p, \delta}$$
  
  Write $f(R, t)=f_R(t).$
  
  Then $$|f_{R_1}(t)-f_{R_2}(t)|=|R_1-R_2|\cdot |\int_{[0,1]}\partial_Rf ((1-s)R_1+sR_2, t)ds|$$ $$ \leq |R_1-R_2|\cdot \int_{[0,1]}|\partial_Rf ((1-s)R_1+sR_2, t)|ds$$ $$ \leq |R_1-R_2|\cdot|\partial_Rf|_{C^0}=C_0|R_1-R_2|$$ so that $\|\tau_{R_1}f_{R_1}-\tau_{R_2}f_{R_2}\|_{C_0}\leq C_0|R_1-R_2|$.

   Hence $$\|[\tau_{R_1}f_{R_1}{\hat \xi}]\circ \tau_{-2R_2}-[\tau_{R_2}f_{R_2}{\hat \xi}]\circ \tau_{-2R_2}\|_{0, p, \delta}\leq  C_0|R_1-R_2| e^{2\delta(|A(R_1)|+|B(R_1)|)} \cdot 
   ||{\hat \xi}\|_{0, p, \delta}.$$

   Therefore $$\|[\tau_{R_1}f_{R_1}{\hat \xi}]\circ
   \tau_{-2R_1}-[\tau_{R_2}f_{R_2}{\hat \xi}]\circ \tau_{-2R_2}\|_{0, p, \delta}\leq C_0|R_1-R_2| e^{2\delta(|A(R_1)|+|B(R_1)|)} \cdot 
   ||{\hat \xi}\|_{0, p, \delta}$$ $$  +|e^{2\delta R_1}-e^{2\delta R_2}|/\delta\cdot C_0\cdot e^{-\delta (R_1)} 
   \|  {\hat \xi}(t)\|_{1, p, \delta}.$$
   
   So that  $$\|H(R_1)-H(R_2)\|_{o}\leq  C_0\{|R_1-R_2| e^{2\delta(|A(R_1)|+|B(R_1)|)}
    +|e^{2\delta R_1}-e^{2\delta R_2}|/\delta\cdot e^{-\delta (R_1)}\} 
   .$$

     \QED

     \begin{cor}
     	$\|D(R)\|_{o}\leq e^{-\delta R}$.

     \end{cor}

     \proof

     $$\|D(R)\|_{o}	\leq \|D(R)-D(R_1)\|_{o}+\|D(R_1)\|_{o}\leq  |e^{-\delta R}-e^{-\delta R_1}|/\delta+\sup_{\|\xi|_{1, p, \delta}\leq 1}
     \|\xi\tau_{R_1}\|_{0, p, \delta}
     $$ $$ \leq |e^{-\delta R}-e^{-\delta R_1}|/\delta+\sup_{\|\xi|_{1, p, \delta}\leq 1}	 e^{-\delta R_1}\|\xi\|_{0, p,\delta}
     =	|e^{-\delta R}-e^{-\delta R_1}|/\delta+ e^{-\delta R_1}.$$

     Now let $R_1\rightarrow \infty$, we get
     $\|D(R)\|_{o}\leq e^{-\delta R}$.

     \QED

     \begin{cor}
     	$D(R)$  extends  continuous over $[R_0, \infty]$ with $D(\infty)=0$.

     \end{cor}

     Back to the function  $F: L_{1, p,\delta}(C_+, E)\times [0, \infty)\rightarrow L_{0, p,\delta}(C_+, E)$ or $F: L_{1, p,\delta}({\bf R}^1\times S^1, E)\times {\bf R}^1 \rightarrow L_{0, p,\delta}({\bf R}^1\times S^1, E)$ defined by $F(u, R)=u\circ \tau_R$. For each fixed $R$, $F$ is linear in $u$ and bounded even in $L_{1, p,\delta}$-norm of the target:
     $F(au_1+bu_2)=(au_1+bu_2)\circ \tau=au_1\circ \tau+ bu_2\circ \tau=aF(u_1)+bF(u_2);$   and $\|F(u)\|_{1, p, \delta}= \|u\circ\tau\|_{1, p, \delta}= e^{\pm\delta(R)}\|u\|_{1, p, \delta}.$
     This proves the following
     \begin{lemma}

     The partial derivative $(D_1F)_{u, R}=:(D_uF)_{u, R}$ is given by   $(D_1F)_{u, R}(\xi)=\xi\circ \tau_R.$
     \end{lemma}

     Next consider  the partial derivative of $F$  along $R$-direction.

     \begin{lemma}
     	The partial derivative $(D_RF)_{u, R}$ is given by   $(D_RF)_{u, R}(\frac{\partial}{\partial R})=\frac{\partial F}{\partial R}={\partial_t u}\circ \tau_R.$
     \end{lemma}

     \proof

     Assume that $u$  is smooth first.
     Then $\|F(u, R+s
     )-F(u, R)-s\cdot {\partial_t u}\circ \tau_R\|_{0, p, \delta}=\|u\circ\tau_{R+s}-u\circ\tau_{R}-s\cdot {\partial_t u}\circ \tau_R\|_{0, p, \delta}=|s|\cdot \|\int _{[0,1]} \{\partial_tu(v(R+s)+(1-v){R})- {\partial_t u}\circ \tau_R\} dv\|_{0, p, \delta}\leq |s|\int_{[0,1]}\| \{\partial_tu(vs+{R})-\cdot {\partial_t u}\circ \tau_R\} \|_{0, p, \delta}dv$.

      We will prove in next lemma that above inequality still true for $u\in L_{1,\delta}^p$.

     In other word, the  inequality before for smooth $u$ still holds for $u\in L_{1,\delta}^p$.

     Then for  $u\in L_{1, p, \delta}$ and $R\in [R_0, \infty)$,  $\lim_{s\rightarrow 0}\int _{[0,1]} \|\{\partial_tu(vs+{R})- {\partial_t u}\circ \tau_R\}\|_{0, p, \delta} dv=\int _{[0,1]}\lim_{s\rightarrow 0} \|\{\partial_tu(vs+{R})- {\partial_t u}\circ \tau_R\} \|_{0, p, \delta}dv=\int _{[0,1]}\|\lim_{s\rightarrow 0} \{\partial_tu_n(vs+{R})- {\partial_t u_n}\circ \tau_R\} \|_{0, p, \delta}dv=0.$ Here the identity before the last one  follows from the dominated convergence theorem.

     \QED

     \begin{lemma}
     The inequality for smooth function $u$ in the proof of the above lemma    still holds for $u\in L_{1,\delta}^p$.
     	
     \end{lemma}

     \proof

The result is well-known. For completeness, we include a proof here.
      For general $u=\lim_{n\rightarrow \infty} u_n$ in $L_{1, \delta}^p$-norm with $u_n$ being smooth,  since the functions $u\rightarrow \|u\|_{1,p, \delta}$ and $(u, R)\rightarrow u\circ\tau_R$ are  continuous,
     $$\|F(u, R+s
     )-F(u, R)-s\cdot {\partial_t u}\circ \tau_R\|_{0, p, \delta}$$ $$ =\|\lim_{n\rightarrow \infty}u_n\circ \tau_{R+s}-\lim_{n\rightarrow \infty}u_n\circ \tau_{R}-s\cdot \lim_{n\rightarrow \infty}{\partial_t u_n}\circ \tau_R\|_{0, p, \delta}$$ $$ =\lim_{n\rightarrow \infty}|s|\cdot \|\int _{[0,1]} \{\partial_tu_n(v(R+s)+(1-v){R})- {\partial_t u_n}\circ \tau_R\} dv\|_{0, p, \delta}$$ $$ \leq |s|\lim_{n\rightarrow \infty}\int_{[0,1]}\| \{\partial_t u_n(vs+{R})- {\partial_t u_n}\circ \tau_R\} \|_{0, p, \delta}dv$$ $$ =|s|\int_{[0,1]}\lim_{n\rightarrow \infty}\| \{\partial_t u_n(vs+{R})- {\partial_t u_n}\circ \tau_R\} \|_{0, p, \delta}dv$$ $$ =|s|\int_{[0,1]}\| \lim_{n\rightarrow \infty}\{\partial_t u_n(vs+{R})- {\partial_t u_n}\circ \tau_R\} \|_{0, p, \delta}dv$$ $$ =|s|\int_{[0,1]}\| \{\partial_tu(vs+{R})-\cdot {\partial_t u}\circ \tau_R\} \|_{0, p, \delta}dv.$$ Here the identity interchanging the integral and limit follows from the theorem  on dominated convergence. Indeed, the function $f_n(v)=:\| \{\partial_t u_n(vs+{R})- {\partial_t u_n}\circ \tau_R\} \|_{0, p, \delta}$ is continuous when $u_n$ is smooth and $|f_n(v)|=f_n(v)\leq  e^{\delta(R+|s|)}\{\| u\|_{1, p, \delta}+1\}$ so that      $\lim_{n\rightarrow \infty}\int_{[0,1]}\| \{\partial_t u_n(vs+{R})- {\partial_t u_n}\circ \tau_R\} \|_{0, p, \delta}dv=\int_{[0,1]}\lim_{n\rightarrow \infty}\| \{\partial_t u_n(vs+{R})- {\partial_t u_n}\circ \tau_R\} \|_{0, p, \delta}dv.$

       \QED

       By the second lemma of this section, $\|(D_RF)_{u, R}(\frac{\partial}{\partial R})\|_{0, p, \delta }=\|{\partial_t u}\circ \tau_R\|_{0, p, \delta }\leq e^{-\delta R}\|u\|_{1, p,\delta}.$
       This implies that $D_RF$ can be extended continuously over $R=\infty$ with  the value equal to zero  over  $R=\infty$.

       Similarly, $(D_\theta F)_{u, R_{\theta}}(\frac{\partial}{\partial \theta})=\partial_s u\circ \tau_{R_{\theta}}$. Then

       $$\|1/r\cdot(D_\theta F)_{u, R_{\theta}}(\frac{\partial}{\partial \theta})||_{k-1, p, \delta}=1/r\cdot \|\partial_s u\circ \tau_{R_{\theta}}||_{k-1, p, \delta}=1/r \cdot \|\partial_s u\circ \tau_{R}||_{k-1, p, \delta}$$ $$ \sim 1/r \cdot e^{-\delta R}\|u\|_{k, p, \delta}\sim 1 \cdot e^{-\delta R/2}\|u\|_{k, p, \delta}.$$
      Combing all the lemmas above, this proves the following

      \begin{cor}
       The function F is of class $C^1$ for  $R\not = \infty$ and $DF$ can be
       extended continuously over $R=\infty.$
      \end{cor}

      In  fact, the formula for $\partial_{u}F$ is  already  valid even for $R=\infty.$

      To prove  $\partial_{x}\Psi_L$ and $\partial_{y}\Psi_L$ in last section or $\partial_{x}\Psi_N$ and $\partial_{y}\Psi_N$
       in the sequel of this paper   are  the real derivatives, we need  the following lemma.  %Note that the worst  scenario in this paper is  $\|\partial_{r}\Psi_L(u, r)\|_{0, p,\delta}\sim 1/\ln^2R \|u\|_{1, p, \delta}\sim f(r) \|u\|_{1, p, \delta}$ with $f(r)>0$  for $r\in (0, r_0)$ with $\lim_{r\rightarrow 0}f(r)=0.$
       Note  that both $\Psi_L$  and $\Psi_N$ are already continuous including $r=0$.

   \begin{lemma}
   Let $K: W\times I =L^k_{ p, \delta}\times (-x_0, x_0)\rightarrow L=L_{k-1, p, \delta} $ be a  continuous function  such that $\partial_x K: W\times (I\setminus \{0\})\rightarrow L $  is continuous and extended continuously over $x=0$, with the extension denoted by  ${{\widetilde \partial_x} K}:W\times I\rightarrow L $. Then $\partial_x K$ exits over  $W\times \{x=0\}$.

   \end{lemma}

    \proof

    For $s>0,$
    $$\|K(u, s)-K(u, 0)-{{\widetilde \partial_x} K}_{(u, 0)}s\|_{k-1, p, \delta}=||K(u, s)-\lim_{\mu>0, \mu\rightarrow 0}K(u, \mu)-{{\widetilde \partial_x} K}_{(u, 0)}s\|_{k-1, p, \delta}$$ $$=s\cdot || \lim_{\mu>0, \mu\rightarrow 0}\int_\mu^1\{\partial_x K(u, vs+(1-v)\mu)-{\widetilde \partial_x} K(u, 0)\}dv\|_{k-1, p, \delta}$$ $$ \leq s\cdot  \lim_{\mu>0, \mu\rightarrow 0}\int_\mu^1||{\widetilde \partial_x} K(u, vs+(1-v)\mu)-{\widetilde \partial_x} K(u, 0)\|_{k-1, p, \delta}dv.$$

    For any given $\epsilon >0, $ the continuity of ${\widetilde \partial_x} K$ implies that there exists a $\rho>0$ such that when $|vs+(1-v)\mu|<\rho$, $$||{\widetilde \partial_x} K(u, vs+(1-v)\mu)-{\widetilde \partial_x} K(u, 0)\|_{k-1, p, \delta}<\epsilon $$ for fixed $u$.  Now for fixed $s>0$, if $0<\mu<s$,    the condition that $s<\rho$ implies that $|vs+(1-v)\mu|<\rho$ so that $$ \lim_{\mu>0, \mu\rightarrow 0}\int_\mu^1||{\widetilde \partial_x} K(u, vs+(1-v)\mu)-{\widetilde \partial_x} K(u, 0)\|_{k-1, p, \delta}dv$$ $$ = \lim_{0<\mu<s, \mu\rightarrow 0}\int_\mu^1||{\widetilde \partial_x} K(u, vs+(1-v)\mu)-{\widetilde \partial_x} K(u, 0)\|_{k-1, p, \delta}dv $$ $$ \leq \epsilon.$$ In other words, the function $E(u, s)=: \lim_{\mu>0, \mu\rightarrow 0}\int_\mu^1||{\widetilde \partial_x} K(u, vs+(1-v)\mu)-{\widetilde \partial_x} K(u, 0)\|_{k-1, p, \delta}dv$ has the property that $\lim_{s\rightarrow 0}E(u, s)=0.$ Similar result holds for $s<0$ so that $\|K(u, s)-K(u, 0)-{{\widetilde \partial_x} K}_{(u, 0)}s\|_{k-1, p, \delta}\sim |s| o(|s|).$

    \QED

    Combining all the results so far, we have proved the main theorem of this paper.

%{\small \smallskip\noindent Updated 5 December 2006.}

\end{document}